\documentclass[a4paper]{article}
\usepackage{amsmath,amsthm,amssymb,amscd,mathrsfs}

\numberwithin{equation}{section}

{\theoremstyle{definition}\newtheorem{definition}{Definition}[section]
\newtheorem{notation}[definition]{Notation}
\newtheorem{terminology}[definition]{Terminology}
\newtheorem{remark}[definition]{Remark}
}
\newtheorem{proposition}[definition]{Proposition}
\newtheorem{lemma}[definition]{Lemma}
\newtheorem{theorem}[definition]{Theorem}
\newtheorem{corollary}[definition]{Corollary}

\newcommand{\Ptil}{\widetilde{P}}
\newcommand{\cK}{\mathcal{K}}
\newcommand{\Aut}{\operatorname{Aut}}

\newcommand{\Out}{\operatorname{Out}}
\newcommand{\si}{\sigma}
\newcommand{\om}{\omega}
\newcommand{\recht}{\rightarrow}
\newcommand{\R}{\mathbb{R}}
\newcommand{\B}{\operatorname{B}}
\newcommand{\M}{\operatorname{M}}
\newcommand{\Mor}{\operatorname{Mor}}
\newcommand{\eps}{\varepsilon}
\newcommand{\Tr}{\operatorname{Tr}}
\newcommand{\GL}{\operatorname{GL}}

\newcommand{\de}{\Delta}
\newcommand{\N}{\mathbb{N}}
\newcommand{\C}{\mathbb{C}}

\newcommand{\be}{\beta}
\newcommand{\ot}{\otimes}
\newcommand{\id}{\operatorname{id}}

\newcommand{\dpr}{^{\prime\prime}}
\newcommand{\la}{\langle}
\newcommand{\ra}{\rangle}
\newcommand{\Sd}{\operatorname{Sd}}
\newcommand{\Ad}{\operatorname{Ad}}
\newcommand{\deh}{\hat{\de}}
\newcommand{\fancyot}{\mathbin{\text{\footnotesize\textcircled{\tiny \sf T}}}}

\newcommand{\T}{\mathbb{T}}
\newcommand{\cst}{C$^*$}

\newcommand{\cL}{\mathcal{L}}

\newcommand{\ltimesfull}{\;_{\text{\rm f}}\hspace{-.6ex}\ltimes}
\newcommand{\ltimesred}{\;_{\text{\rm r}}\hspace{-.6ex}\ltimes}
\newcommand{\epsh}{\hat{\eps}}

\newcommand{\otmin}{\underset{\text{\rm min}}{\otimes}}
\newcommand{\otmax}{\underset{\text{\rm max}}{\otimes}}

\newcommand{\SU}{\operatorname{SU}}
\newcommand{\U}{\operatorname{U}}
\newcommand{\cG}{\mathbb{G}}
\newcommand{\cGh}{\widehat{\mathbb{G}}}
\newcommand{\vphi}{\varphi}
\newcommand{\cB}{\mathcal{B}}
\newcommand{\cBtil}{\widetilde{\mathcal{B}}}
\newcommand{\Btil}{\widetilde{B}}
\newcommand{\omtil}{\widetilde{\omega}}
\newcommand{\qbinom}[2]{\bigl[\begin{smallmatrix}\displaystyle #1 \\
    \displaystyle #2 \end{smallmatrix}\bigr]}
\newcommand{\slim}{\operatornamewithlimits{s\mbox{$^*$}-lim}}
\newcommand{\supp}{\operatorname{supp}}
\newcommand{\Irred}{\operatorname{Irred}}
\newcommand{\dimq}{\operatorname{dim}_q}
\newcommand{\ox}{\overline{x}}
\newcommand{\ovt}{\overline{\otimes}}
\newcommand{\red}{_{\text{\rm red}}}
\newcommand{\vnalg}{_{\text{\rm red}}^{\prime\prime}}
\newcommand{\hhl}{\widehat{h}_L}
\newcommand{\hhr}{\widehat{h}_R}
\newcommand{\Lambdahl}{\widehat{\Lambda}_L}
\newcommand{\Lambdahr}{\widehat{\Lambda}_R}
\newcommand{\End}{\operatorname{End}}
\newcommand{\dehop}{\deh^{\text{\rm op}}}

\newcommand{\Vb}{\mathbb{V}}  %The universal unitary between discrete
\newcommand{\cV}{\mathscr{V}} %Right regular of compact
\newcommand{\cW}{\mathscr{W}} %Left regular of compact
\newcommand{\cVh}{\widehat{\mathscr{V}}} %Right regular of discrete
\newcommand{\cWh}{\widehat{\mathscr{W}}} %Left regular of discrete

\newcommand{\lambdah}{\widehat{\lambda}}
\newcommand{\rhoh}{\widehat{\rho}}

\newcommand{\rhoop}{\rho^\text{\rm op}}
\newcommand{\phiml}{{\phi_l^-}}
\newcommand{\phimr}{{\phi_r^-}}
\newcommand{\phipl}{{\phi_l^+}}
\newcommand{\phipr}{{\phi_r^+}}

\newcommand{\Calg}{C_{\text{\rm alg}}}

\newcommand{\Z}{\mathbb{Z}}
\newcommand{\Inn}{\operatorname{Inn}}

\newcommand{\xibol}{\xi\hspace{-1ex}\raisebox{1.6ex}[0pt][0pt]{$\scriptstyle\circ$}}
\newcommand{\otalg}{\otimes_{\text{\rm alg}}}

\newcommand{\deltah}{\widehat{\delta}}
\newcommand{\omu}{\overline{\mu}}
\newcommand{\Atil}{\widetilde{A}}

\begin{document}
\begin{center}
{\LARGE\bf The boundary of universal discrete quantum groups, exactness and factoriality}

\renewcommand{\thefootnote}{}
\footnotetext{This research was supported by the EU-Network
Quantum-Spaces --- Noncommutative Geometry (HPRN-CT-2002-00280).}

\bigskip

{\sc by Stefaan Vaes$^{\text{\rm (a,b)}}$ and Roland Vergnioux $^{\text{\rm (c,d)}}$}
\end{center}

{\footnotesize \parbox[t]{5mm}{(a)}CNRS; Institut de Math{\'e}matiques de Jussieu; 175, rue du Chevaleret; F--75013 Paris (France) \\
\parbox[t]{5mm}{(b)}K.U.Leuven; Department of Mathematics; Celestijnenlaan 200B;
B--3001 Leuven (Belgium) \\
\parbox[t]{5mm}{(c)}Universit{\"a}t M{\"u}nster; Institute of Mathematics; Einsteinstra{\ss}e 62;
D--48149 M{\"u}nster (Germany) \\
\parbox[t]{5mm}{(d)}Universit{\'e} de Caen; Laboratoire de Math{\'e}matiques Nicolas Oresme; BP 5186;
F--14032 Caen Cedex (France)

e-mail: stefaan.vaes@wis.kuleuven.be, roland.vergnioux@math.unicaen.fr

}

\begin{abstract}
\noindent We study the \cst-algebras and von Neumann algebras
associated with the universal discrete quantum groups. They give rise
to full prime factors and simple exact \cst-algebras. The main tool in
our work is the study of an amenable boundary action, yielding the
Akemann-Ostrand property. Finally, this boundary can be identified
with the Martin or the Poisson boundary of a quantum random walk.
\end{abstract}

\section*{Introduction}

Since Murray and von Neumann introduced \emph{von Neumann algebras},
the ones associated with discrete groups played a prominent role. The
main aim of this paper
is to show how concrete examples of discrete quantum groups give rise
to interesting \cst- and von Neumann algebras.

In the 1980's, Woronowicz \cite{wor2} introduced the notion of a \emph{compact quantum group} and generalized the classical Peter-Weyl representation theory.
Many fascinating examples of compact quantum groups are available by now: Drinfel'd and Jimbo \cite{drinfeld,jimbo} introduced the
\emph{$q$-deformations of compact semi-simple Lie groups} and Rosso \cite{rosso} showed that they fit into the theory of Woronowicz. The
\emph{universal orthogonal and unitary} quantum groups were introduced by Van Daele and Wang \cite{VDW} and studied in detail by Banica
\cite{banica1,banica2}.

This paper mainly deals with the universal orthogonal quantum group $\cG = A_o(F)$, defined from a matrix $F \in \GL(n,\C)$ satisfying $F
\overline{F} = \pm 1$. Its underlying \cst-algebra $C(\cG)$ is the universal \cst-algebra generated by the entries of a unitary $n$ by $n$ matrix
$(U_{ij})$ satisfying $(U_{ij}) = F(U_{ij}^*) F^{-1}$. Using the GNS construction of the (unique) Haar state of $\cG$, we obtain the reduced
\cst-algebra $C(\cG)\red$ and the von Neumann algebra $C(\cG)\vnalg$. This paper deals with a detailed study of these operator algebras. Note that,
for $n \geq 3$, $C(\cG)\red$ is a non-trivial quotient of $C(\cG)$ by non-amenability of the discrete quantum group $\cGh$.

In Section \ref{sec.boundary}, we construct a \emph{boundary} for the
dual $\cGh$ of $\cG= A_o(F)$. This boundary $\cB_\infty$ is a unital
\cst-algebra that admits a natural action of $\cGh$. In Section
\ref{sec.amenable}, we introduce the notion of an \emph{amenable
  action} of a discrete quantum group on a unital \cst-algebra. This
definition involves a non-trivial algebraic condition, which is the
proper generalization of Anantharaman-Delaroche's centrality condition
(see Th{\'e}or{\`e}me 3.3 in \cite{AD}). We then prove that the boundary
action of $\cGh$ is amenable. The construction of the boundary
$\cB_\infty$ and the proof of the amenability of the boundary action
involve precise estimates on the representation theory of
$A_o(F)$. These estimates are dealt with in the appendix.

From the amenability of the boundary action of the dual of $\cG = A_o(F)$, we deduce that the reduced \cst-algebra $C(\cG)\red$ is \emph{exact} and
satisfies the \emph{Akemann-Ostrand property}. In the setting of finite von Neumann algebras, Ozawa \cite{ozawa} showed that the Akemann-Ostrand
property implies solidity of the associated von Neumann algebra. Since in general $C(\cG)\vnalg$ is of type III, we need a generalization of Ozawa's
definition (see Section \ref{sec.solid}) and we deduce that, for $\cG= A_o(F)$, the von Neumann algebras $C(\cG)\vnalg$ are \emph{generalized solid}.
In particular, for $F$ the $n$ by $n$ identity matrix, we get a solid von Neumann algebra.

In Section \ref{sec.probabilistic}, we make the link between our boundary $\cB_\infty$ for the dual of $\cG = A_o(F)$ and boundaries arising from
quantum random walks on $\cGh$. We construct a \emph{harmonic state} $\om_\infty$ on $\cB_\infty$ and identify $(\cB_\infty,\om_\infty)\dpr$ with the
\emph{Poisson boundary} of a random walk on $\cGh$. Poisson boundaries of discrete quantum groups were defined by Izumi in \cite{I}, who computed
them for the dual of $\SU_q(2)$. This computation was then extended to the dual of $\SU_q(n)$ in \cite{INT}. We identify $\cB_\infty$ with the
\emph{Martin boundary} of a random walk on $\cGh$. Martin boundaries of discrete quantum groups were defined by Neshveyev and Tuset in \cite{NT} and
computed there for the dual of $\SU_q(2)$.

In the short Section \ref{sec.generalexact}, we provide a general
exactness result for quantum group \cst-algebras $C(\cG)\red$. We show
that for \emph{monoidally equivalent} quantum groups $\cG$ and $\cG_1$
(see \cite{BDRV}), $C(\cG)\red$ is exact if and only if $C(\cG_1)\red$
is exact. This provides an alternative proof for the exactness of
$C(\cG)\red$ when $\cG = A_o(F)$ and proves exactness in other
examples as well.

In Section \ref{sec.factoriality}, we deal with \emph{factoriality} of the von Neumann algebra $C(\cG)\vnalg$ and \emph{simplicity} of the
\cst-algebra $C(\cG)\red$, whenever $\cG = A_o(F)$ and $F$ is at least a $3$ by $3$ matrix. We were only able to settle factoriality and simplicity
assuming an extra condition on the norm of $F$: if $\sqrt{5} \|F\|^2 \leq \Tr(F^*F)$, the von Neumann algebra $C(\cG)\vnalg$ is a \emph{full factor}
and we compute its Connes invariants. If $8 \|F\|^8 \leq 3 \Tr(F^*F)$, the \cst-algebra $C(\cG)\red$ is simple. Both conditions are satisfied when
$F$ is sufficiently close to the $n$ by $n$ identity matrix, for $n \geq 3$. Moreover, it is our belief that they are superfluous. Note that
simplicity of the reduced \cst-algebra of the universal unitary quantum groups $A_u(F)$ was proven by Banica in \cite{banica2}. For $\cG = A_u(F)$,
the fusion algebra can be described using the free monoid $\N \ast \N$, while for $\cG = A_o(F)$, the fusion algebra is the same as the one of
$\SU(2)$ and is, in particular, abelian. For that reason, Banica's approach is closer to Powers proof of the simplicity of $C^*_r(\mathbb{F}_n)$.

For the convenience of the reader, we included a rather extensive
section of preliminaries, dealing with the general theory of compact/discrete
quantum groups, their actions on \cst-algebras and exactness.

\section{Preliminaries}

We use the symbol $\ot$ to denote several types of tensor products. In particular $\ot$ denotes the \emph{minimal} tensor product of \cst-algebras.
The \emph{maximal} tensor product is denoted by $\otmax$. If we want to stress the difference with the minimal tensor product, we write $\otmin$. We
also make use of the leg numbering notation in multiple tensor products: if $a \in A \ot A$, then $a_{12},a_{13},a_{23}$ denote the obvious elements
in $A \ot A \ot A$, e.g.\ $a_{12} = a \ot 1$.

\subsection*{Compact quantum groups}

We briefly overview the theory of compact quantum groups developed by
Woronowicz in \cite{wor2}. We refer to the survey paper \cite{MVD} for
a smooth approach to these results.

\begin{definition}
A \emph{compact quantum group} $\cG$ is a pair $(C(\cG),\de)$, where
\begin{itemize}
\item $C(\cG)$ is a unital \cst-algebra;
\item $\de : C(\cG) \recht C(\cG) \ot C(\cG)$ is a unital
  $^*$-homomorphism satisfying the \emph{co-associativity} relation
$$(\de \ot \id)\de = (\id \ot \de)\de \; ;$$
\item $\cG$ satisfies the \emph{left and right cancellation property}
  expressed by
$$\de(C(\cG))(1 \ot C(\cG)) \quad\text{and}\quad \de(C(\cG))(C(\cG) \ot 1)
\quad\text{are total in}\;\; C(\cG) \ot C(\cG) \; .$$
\end{itemize}
\end{definition}

\begin{remark}
We use the fancy notation $C(\cG)$ to suggest the analogy with the
basic example given by continuous functions on a compact group. In the
quantum case however, there is no underlying space $\cG$ and $C(\cG)$
is a non-abelian \cst-algebra.
\end{remark}

The two major aspects of the general theory of compact quantum groups
are the existence and uniqueness of a Haar measure and the Peter-Weyl
representation theory.

\begin{theorem}
Let $\cG$ be a compact quantum group. There exists a unique state $h$ on $C(\cG)$ satisfying $(\id \ot h)\de(a) = h(a)1 = (h \ot \id)\de(a)$ for all
$a \in C(\cG)$. The state $h$ is called the \emph{Haar state} of $\cG$.
\end{theorem}

\begin{definition}
A \emph{unitary representation} $U$ of a compact quantum group $\cG$
on a Hilbert space $H$ is a unitary element $U \in \M(\cK(H) \ot
C(\cG))$ satisfying
\begin{equation} \label{eq.rep}
(\id \ot \de)(U) = U_{12} U_{13} \; .
\end{equation}
Whenever $U^1$ and $U^2$ are unitary representations of $\cG$ on the
respective Hilbert spaces $H_1$ and $H_2$, we
define
$$\Mor(U^1,U^2) := \{ T \in \B(H_2,H_1) \mid U_1(T \ot 1) = (T \ot
1)U_2 \}\; .$$ The elements of $\Mor(U^1,U^2)$ are called \emph{intertwiners}. We use the notation $\End(U) := \Mor(U,U)$. A unitary representation
$U$ is said to be \emph{irreducible} if $\End(U) = \C1$. Two unitary representations $U^1$ and $U^2$ are said to be \emph{unitarily equivalent} when
$\Mor(U^1,U^2)$ contains a unitary operator.
\end{definition}

The following result is crucial.

\begin{theorem}
Every irreducible representation of a compact quantum group is
finite-dimensional. Every unitary representation is unitarily
equivalent to a direct sum of irreducibles.
\end{theorem}

Because of this theorem, we almost exclusively deal with
finite-dimensional representations, the regular representation being
the exception. By choosing an orthonormal basis
of the Hilbert space $H$, a finite-dimensional unitary representation of $\cG$ can be
considered as a unitary matrix $(U_{ij})$ with entries in $C(\cG)$ and
\eqref{eq.rep} becomes
$$\de(U_{ij}) = \sum_k U_{ik} \ot U_{kj} \; .$$

The product in the \cst-algebra $C(\cG)$ yields a tensor product on
the level of unitary representations.

\begin{definition}
Let $U^1$ and $U^2$ be unitary representations of $\cG$ on the
respective Hilbert spaces $H_1$ and $H_2$. We define the tensor
product
$$U^1 \fancyot U^2 := U^1_{13} U^2_{23} \in \M(\cK(H_1 \ot H_2) \ot
C(\cG)) \; .$$
\end{definition}

\begin{notation}
We denote by $\Irred(\cG)$ the set of (equivalence classes) of
irreducible unitary representations of a compact quantum group
$\cG$. We choose representatives $U^x$ on the Hilbert space $H_x$ for
every $x \in \Irred(\cG)$. Whenever $x,y \in \Irred(\cG)$, we use $x
\ot y$ to denote the unitary representation $U^x \fancyot U^y$. The
class of the trivial unitary representation is denoted by $\eps$.
\end{notation}

The set $\Irred(\cG)$ is equipped with a natural involution $x \mapsto
\ox$ such that $U^{\ox}$ is the unique (up to unitary equivalence)
irreducible unitary representations such that
$$\Mor(x \ot \ox,\eps) \neq 0 \neq \Mor(\ox \ot x,\eps) \; .$$
The unitary representation $U^{\ox}$ is called the contragredient of
$U^x$.

The irreducible representations of $\cG$ and the Haar state $h$ are
connected by the \emph{orthogonality relations}. For every $x \in
\Irred(\cG)$, we have a unique invertible positive self-adjoint
element $Q_x \in \B(H_x)$ satisfying $\Tr(Q_x) = \Tr(Q_x^{-1})$ and
\begin{equation}\label{eq.orthogonality}
(\id \ot h)(U^x (A \ot 1) (U^x)^*) = \frac{\Tr(Q_xA)}{\Tr(Q_x)} 1 \quad , \quad
(\id \ot h)((U^x)^*(A \ot 1) U^x) = \frac{\Tr(Q_x^{-1}A)}{\Tr(Q_x^{-1})} 1
\; ,
\end{equation}
for all $A \in \B(H_x)$.

\begin{definition}
For $x \in \Irred(\cG)$, the value $\Tr(Q_x)$ is called the
\emph{quantum dimension} of $x$ and denoted by $\dimq(x)$. Note that
$\dimq(x) \geq \dim(x)$, with equality holding if and only if $Q_x = 1$.
\end{definition}

\subsection*{Discrete quantum groups and duality}

A discrete quantum group is defined as the dual of a compact quantum
group, putting together all irreducible representations.

\begin{definition}
Let $\cG$ be a compact quantum group. We define the dual (discrete)
quantum group $\cGh$ as follows.
\begin{align*}
c_0(\cGh) &= \bigoplus_{x \in \Irred(\cG)} \B(H_x) \; , \\
\ell^\infty(\cGh) &= \prod_{x \in \Irred(\cG)} \B(H_x)
\end{align*}
We denote the minimal central projections of $\ell^\infty(\cGh)$ by
$p_x$, $x \in \Irred(\cG)$.

We have a natural unitary $\Vb \in \M(c_0(\cGh) \ot C(\cG))$ given by
$$\Vb = \bigoplus_{x \in \Irred(\cG)} U^x \; .$$
The unitary $\Vb$ implements the duality between $\cG$ and $\cGh$.
We have a natural comultiplication
$$\deh : \ell^\infty(\cGh) \recht
\ell^\infty(\cGh) \ovt \ell^\infty(\cGh) : (\deh \ot \id)(\Vb) = \Vb_{13} \Vb_{23} \; .$$\end{definition}

The notation introduced above is aimed to suggest the basic example
where $\cG$ is the dual of a discrete group $\Gamma$, given by $C(\cG)
= C^*(\Gamma)$ and $\de(\lambda_x) = \lambda_x \ot \lambda_x$ for all
$x \in \Gamma$. The map $x \mapsto \lambda_x$ yields an identification
of $\Gamma$ and $\Irred(\cG)$ and then, $\ell^\infty(\cGh) = \ell^\infty(\Gamma)$.

\begin{remark}
It is of course possible to give an intrinsic definition of a discrete
quantum group (not as the dual of a compact quantum group). This was
already implicitly clear in Woronowicz' work and was explicitly done
in \cite{ER,VD}. For our purposes, it is most convenient to take the
compact quantum group as a starting point: indeed, all interesting
examples of concrete discrete quantum groups are defined as the dual
of certain compact quantum groups. We shall study one particular class below.
\end{remark}

The discrete quantum group $\ell^\infty(\cGh)$ comes equipped with a natural modular structure.

\begin{notation} \label{not.states}
We choose \emph{unit} vectors $t_x \in H_x \ot H_{\ox}$ invariant
under $U^x \fancyot U^{\ox}$. The vectors $t_x$ are unique up to
multiplication by $\T$. We then have canonically defined states
$\vphi_x$ and $\psi_x$ on $\B(H_x)$ related to
\eqref{eq.orthogonality} as follows.
\begin{align*}
\psi_x(A) &= t_x^* (A \ot 1) t_x = \frac{\Tr(Q_xA)}{\Tr(Q_x)} = (\id \ot h)(U^x (A \ot 1) (U^x)^*) \quad\text{and} \\ \vphi_x(A) &= t_{\ox}^* (1 \ot
A) t_{\ox} = \frac{\Tr(Q_x^{-1}A)}{\Tr(Q_x^{-1})} = (\id \ot h)((U^x)^*(A \ot 1) U^x) \; ,
\end{align*}
for all $A \in \B(H_x)$. As a complement to the vectors $t_x$, we also choose \emph{unit}
vectors $s_x \in \Mor(x \ot \ox,\eps)$ normalized such that $(s_x^*
\ot 1)(1 \ot t_{\ox}) = \frac{1}{\dimq(x)}$ for all $x \in \Irred(\cG)$.
In certain examples, one can consistently choose $s_x = t_x$, but
this is not always the case.
\end{notation}

The states $\vphi_x$ and $\psi_x$ are significant, since they provide
a formula for the invariant weights on $\ell^\infty(\cGh)$.

\begin{proposition} \label{prop.invariantweights}
The left invariant weight $\hhl$ and the right invariant
weight $\hhr$ on $\cGh$ are given by
$$\hhl = \sum_{x \in \Irred(\cG)} \dimq(x)^2 \psi_x \quad\text{and}\quad  \hhr = \sum_{x \in \Irred(\cG)} \dimq(x)^2 \vphi_x
 \; .$$
\end{proposition}

The following formula is used several times in the paper.

\begin{proposition}
Let $x,y \in \Irred(\cG)$ and suppose that $p^{x \ot y}_z \in \End(x
\ot y)$ is an orthogonal projection onto a subrepresentation
equivalent with $z \in \Irred(\cG)$. Then,
$$(\id \ot \psi_y)(p^{x \ot y}_z) = \frac{\dimq(z)}{\dimq(x)\dimq(y)}
1 \quad\text{and}\quad (\vphi_x \ot \id)(p^{x \ot y}_z) = \frac{\dimq(z)}{\dimq(x)\dimq(y)}
1 \; .$$
\end{proposition}
\begin{proof}
Since $(\id \ot \psi_y)(p^{x \ot y}_z) = (1 \ot t_y^*)(p^{x \ot y}_z \ot 1)(1 \ot t_y) \in \End(x) = \C1$, it suffices to check that \linebreak
$(\psi_x \ot \psi_y)(p^{x \ot y}_z) = \frac{\dimq(z)}{\dimq(x)\dimq(y)}$, which immediately follows from the formula $(Q_x \ot Q_y)T = T Q_z$ for all
$T \in \Mor(x \ot y,z)$.
\end{proof}

\subsection*{Regular representations}

Both the algebras $C(\cG)$ and $c_0(\cGh)$ have two natural
representations on the same Hilbert space.

Using \eqref{eq.orthogonality}, we canonically identify the GNS
Hilbert space $L^2(C(\cG),h)$ with
$$L^2(\cG) := \bigoplus_{x \in \Irred(\cG)} (H_x \ot H_{\ox})$$
by taking
\begin{align}
\rho & : C(\cG) \recht \B(L^2(\cG)) : \rho\bigl( (\om_{\eta,\xi} \ot
\id)(U^x) \bigr) \xi_0 = \xi \ot (1 \ot \eta^*)t_{\ox} \; , \label{eq.regrho} \\
\lambda & : C(\cG) \recht \B(L^2(\cG)) : \lambda\bigl( (\om_{\eta,\xi}
\ot \id)(U^x) \bigr) \xi_0 = (1 \ot \eta^*)t_{\ox} \ot \xi \notag
\end{align}
for all $x \in \Irred(\cG)$ and all $\xi,\eta
\in H_x$. Here $\xi_0$ denotes the canonical unit vector in $H_\eps
\ot H_\eps = \C$. We use the notation $\om_{\eta,\xi}(a) = \langle
\eta,a \xi \rangle$ and we use inner products that are linear in the
second variable.

\begin{notation}
Let $\cG$ be a compact quantum group with Haar state $h$. We denote by
$C(\cG)\red$ the \cst-algebra $\rho(C(\cG))$ given by the GNS
construction for $h$. We denote by $C(\cG)\vnalg$ the generated von
Neumann algebra.
\end{notation}

The aim of this paper is a careful study of $C(\cG)\red$ and
$C(\cG)\vnalg$ for certain concrete examples of compact quantum
groups. The previous paragraph clearly suggest that these algebras are
the natural counterparts of $C^*_r(\Gamma)$ and $\cL(\Gamma)$ for a
discrete group $\Gamma$.

We introduce the left and right regular representations for $\cG$ and
$\cGh$. For the convenience of the reader, we provide several explicit
formulae.

\begin{definition}
The right regular representation $\cV \in \cL(L^2(\cG) \ot C(\cG))$ of
$\cG$ and the left regular representation $\cW \in \cL(C(\cG) \ot
L^2(\cG))$ are defined as
\begin{align*}
\cV (\rho(a) \xi_0 \ot 1) &= \bigl((\rho \ot \id)\de(a)\bigr) (\xi_0 \ot 1) \\
\cW^*(1 \ot \rho(a)\xi_0) &= \bigl((\id \ot \rho)\de(a)\bigr) (1 \ot \xi_0) \; .
\end{align*}
\end{definition}

Recall that $\Vb \in \M(c_0(\cGh) \ot C(\cG))$ is given by $\displaystyle\Vb =
\bigoplus_{x \in \Irred(\cG)} U^x$.

\begin{notation}
We define
\begin{align*}
\lambdah &: \ell^\infty(\cGh) \recht \B(L^2(\cG)) : \lambdah(a) \xi_x =
(a_x \ot 1)\xi_x \quad\text{for all}\;\; a \in \ell^\infty(\cGh),
\xi_x \in H_x \ot H_{\ox} \; , \\
\rhoh &: \ell^\infty(\cGh) \recht \B(L^2(\cG)) : \rhoh(a) \xi_x =
(1 \ot a_{\ox})\xi_x \quad\text{for all}\;\; a \in \ell^\infty(\cGh),
\xi_x \in H_x \ot H_{\ox} \; .
\end{align*}
We define the unitary $u \in \B(L^2(\cG))$ by $u(\xi \ot \eta) = \eta
\ot \xi$ for $\xi \in H_x$, $\eta \in H_{\ox}$. Note that $\rhoh = (\Ad
u) \lambdah$ and that $u^2 = 1$.
\end{notation}

\begin{proposition}
The left and the right regular representation of $\cG$ are given by
$$\cV = (\lambdah \ot \id)(\Vb)  \quad\text{and}\quad \cW = (\id \ot
\rhoh)(\Vb_{21}) \; .$$
\end{proposition}

So, as it should be, the left and the right regular representation of $\cG$ give rise to two commuting representations of $\ell^\infty(\cGh)$. We now
symmetrically and explicitly write down how the left and the right regular representation of $\cGh$ give rise to two commuting representations of
$C(\cG)$, see Proposition \ref{prop.formularegular}.

We explicitly perform the GNS construction for the weights $\hhl$ and
$\hhr$ (see Proposition \ref{prop.invariantweights}) in order to give formulae for the left and right regular
representation of $\cGh$. Recall the choice of unit vectors $s_x$ made
in Notation \ref{not.states}.

\begin{notation} \label{not.GNSdiscrete}
Let $a \in \ell^\infty(\cGh)$. We define, whenever the right hand side
makes sense,
$$\Lambdahl(a) = \sum_{x \in \Irred(\cG)} \dimq(x) (ap_x \ot 1)s_x
\quad\text{and}\quad
\Lambdahr(a) = \sum_{x\in \Irred(\cG)} \dimq(x) u (1 \ot ap_x)s_{\ox}
\; .$$
The maps $\Lambdahl$ and $\Lambdahr$, together with the representation
$\lambdah : \ell^\infty(\cGh) \recht \B(L^2(\cG))$ provide a GNS
construction for $\hhl$ and $\hhr$ respectively.
\end{notation}

\begin{definition}
The left regular representation $\cWh \in \cL(c_0(\cGh) \ot L^2(\cG))$ and the right regular representation $\cVh \in \cL(L^2(\cG) \ot c_0(\cG))$ are
defined as
$$\cWh^*(1 \ot
\Lambdahl(a)) = (\id \ot \Lambdahl)\deh(a) \quad\text{and}\quad \cVh (\Lambdahr(a) \ot 1) = (\Lambdahr \ot \id)\deh(a) \; ,$$ for all $a \in
\ell^\infty(\cGh)$ where $\Lambdahl(a)$, resp.\ $\Lambdahr(a)$ makes sense.
\end{definition}

\begin{proposition} \label{prop.formularegular}
The left and the right regular representation of $\cGh$ are
respectively given by
$$\cWh = (\id \ot \rho)(\Vb) \quad\text{and}\quad \cVh = (\lambda \ot \id)(\Vb_{21})$$
and $\rho = (\Ad u)\lambda$. In particular, the \cst-algebras $\lambda(C(\cG))$ and $\rho(C(\cG))$ commute with each other.
\end{proposition}

\begin{remark}
Our notations and conventions agree with Baaj and Skandalis' \cite{BS}
ones in the following way. We consider $\rho$ as the \lq canonical\rq\
representation of $C(\cG)$ and $\lambdah$ as the \lq canonical\rq\ one
for $c_0(\cGh)$. If we then write $V := (\lambdah \ot \rho)(\Vb) =
(\id \ot \rho)(\cV)$, the
operator $V$ is a \emph{multiplicative unitary} on
$L^2(\cG)$. Together with the unitary $u$, $V$ is \emph{irreducible} in the
sense of D{\'e}finition 6.2 in \cite{BS} and the corresponding
multiplicative unitaries of \cite{BS} are given by
$$\widehat{V} = (\rho \ot \rhoh)(\Vb_{21})  \quad\text{and}\quad
\widetilde{V} = (\lambda \ot \lambdah)(\Vb_{21}) \; .$$
\end{remark}

\subsection*{Actions and crossed products}

We provide a brief introduction to the theory of actions of compact
and discrete quantum groups on \cst-algebras. For details and proofs, see \cite{BS}.

\begin{definition}
A (right) \emph{action of a compact quantum group} $\cG$ on a \cst-algebra $A$ is a non-degenerate $^*$-homomorphism
$$\alpha : A \recht \M(A \ot C(\cG)) \quad\text{satisfying}\quad (\alpha \ot \id)\alpha = (\id \ot \de)\alpha$$
and such that $\alpha(A)(1 \ot C(\cG))$ is total in $A \ot C(\cG)$.

The \emph{crossed product} $A \rtimes \cG$ is defined as the closed linear span of $(\id \ot \rho)\alpha(A) \; (1 \ot \lambdah(c_0(\cG)))$ and is a
\cst-algebra. Observe that $A \rtimes \cG$ is realized as a subalgebra of $\cL(A \ot L^2(\cG))$.
\end{definition}

Note that because of amenability of $\cG$ there is no need to define full and reduced crossed products.

\begin{remark}
The action $\alpha$ of a compact quantum group on a unital \cst-algebra $A$ is said to be \emph{ergodic} if
$$\alpha(a) = a \ot 1 \quad\text{if and only if}\quad a \in \C1 \; .$$
For ergodic actions of compact quantum groups, the usual theory of spectral subspaces is available. In particular, one defines the multiplicity with
which an irreducible representation $x \in \Irred(\cG)$ appears in an ergodic action. We need this notion at one place in the paper and refer to the
introduction of \cite{BDRV} for details.
\end{remark}

\begin{definition}
A (left) \emph{action of a discrete quantum group} $\cGh$ on a \cst-algebra $A$ is a non-degenerate $^*$-homomorphism
$$\alpha : A \recht \M(c_0(\cGh) \ot A) \quad\text{satisfying}\quad (\id \ot \alpha)\alpha = (\deh \ot \id)\alpha$$
and such that $(\epsh \ot \id)\alpha(a) = a$ for all $a \in A$.
\end{definition}

Since $\cGh$ need not be amenable, we introduce the notions of a
covariant representation, full crossed product and reduced crossed
product.

\begin{definition}
Let $\alpha : A \recht \M(c_0(\cGh) \ot A)$ be an action of a discrete
quantum group $\cGh$ on a \cst-algebra $A$.
A \emph{covariant representation} of $(A,\alpha)$ into a
  \cst-algebra $B$ is a pair $(\theta,X)$ where $\theta : A \recht
  \M(B)$ is a non-degenerate $^*$-homomorphism and $X \in \M(c_0(\cGh)
  \ot B)$ is a unitary representation of $\cGh$ satisfying the
  covariance relation
$$(\id \ot \theta)\alpha(a) = X^*(1 \ot \theta(a))X \quad\text{for
  all}\quad a \in A \; .$$
\end{definition}

\begin{proposition}
Let $\alpha : A \recht \M(c_0(\cGh) \ot A)$ be an action of a discrete
quantum group $\cGh$ on a \cst-algebra~$A$.
\begin{itemize}
\item For any covariant representation $(\theta,X)$ of $(A,\alpha)$
  the closed linear span of $$\theta(A) \{(\om \ot \id)(X) \mid \om \in
  \ell^\infty(\cGh)_*\}$$ is a \cst-algebra: the \cst-algebra generated by $(\theta,X)$.
\item The \emph{reduced crossed product} $\cGh \ltimesred A$ is the
  \cst-algebra generated by the \emph{regular covariant
  representation} into $\cL(L^2(\cG) \ot A)$ given by
  $((\lambdah \ot \id)\alpha,\cWh_{12})$.
\item The \emph{full crossed product} $\cGh \ltimesfull A$ is the
  unique (up to isomorphism) \cst-algebra $B$ generated by a covariant
  representation $(\theta,X)$ into $B$ satisfying the following universal
  property: for any covariant representation $(\theta_1,X_1)$ into a
  \cst-algebra $B_1$, there exists a non-degenerate $^*$-homomorphism $\pi : B \recht
  \M(B_1)$ satisfying $\theta_1 = \pi \theta$ and $X_1 = (\id \ot \pi)(X)$.
\end{itemize}
\end{proposition}

\subsection{Exactness}

Recall that a \cst-algebra $A$ is said to be \emph{exact} if the
operation $A \otmin \cdot$ transforms short exact sequences into short
exact sequences.

\begin{definition}
A discrete quantum group $\cGh$ is said to be \emph{exact} if the
operation $\cGh \ltimesred \cdot$ transforms $\cGh$-equivariant short
exact sequences into short exact sequences.
\end{definition}

The following proposition is proved using a classical trick.

\begin{proposition} \label{prop.quantumexact}
A discrete quantum group $\cGh$ is exact if and only if $C(\cG)\red$
is an exact \cst-algebra.
\end{proposition}

\begin{proof}
One implication is obvious by applying the definition of exactness of
$\cGh$ to short exact sequences equivariant with respect to the
trivial action of $\cGh$.

So, suppose that $C(\cG)\red$ is an exact \cst-algebra and suppose that $0 \recht J \recht A \recht A/J \recht 0$ is a $\cGh$-equivariant short exact
sequence. Denote by $\delta_J$, resp.\ $\delta_A$, the actions of $\cGh$ on $J$, resp.\ $A$. For any \cst-algebra $B$ with an action $\delta$ of
$\cGh$ on $B$, we have a canonical injective $^*$-homomorphism
$$\deltah : \cGh \ltimesred B \recht C(\cG)\red \otmin (\cGh
\ltimesfull B) \; ,$$
which is a form of the dual action of $\cG$ on the crossed product.
Observe that at the right hand side, the full crossed product appears.
Consider the commutative diagram
$$
\begin{CD}
0 @>>> \cGh \ltimesred J @>>> \cGh \ltimesred A @>>> \cGh \ltimesred
\frac{A}{J} @>>> 0 \\
@. @VV{\deltah_J}V @VV{\deltah_A}V @VVV @. \\
0 @>>> C(\cG)\red \otmin (\cGh \ltimesfull J) @>>> C(\cG)\red \otmin
(\cGh \ltimesfull A) @>>> C(\cG)\red \otmin (\cGh \ltimesfull \frac{A}{J}) @>>> 0
\end{CD}
$$
By universality, the sequence $0 \recht \cGh \ltimesfull J \recht \cGh
\ltimesfull A \recht \cGh \ltimesfull \frac{A}{J} \recht 0$ is
exact. So, by exactness of $C(\cG)\red$ the bottom row of the
commutative diagram is exact. Suppose $a \in \cGh \ltimesred A$
becomes $0$ in $\cGh \ltimesred \frac{A}{J}$. Since the bottom row of
the diagram is exact, $\deltah_A(a) \in C(\cG)\red \otmin (\cGh
\ltimesfull J)$. Using an approximate identity $(e_\alpha)$ for $J$,
it is easy to check that $\deltah_J(e_\alpha)\deltah_A(a) \recht
\delta_A(a)$. Since $\deltah_A$ is isometric, it follows that $e_\alpha
a \recht a$ and hence, $a \in \cGh \ltimesred J$. This proves the
exactness of the top row in the diagram.
\end{proof}

\subsection*{Examples: the universal compact quantum groups}

The universal compact quantum groups were introduced by Wang and Van
Daele in \cite{VDW}. They are defined as follows.

\begin{definition}
Let $F \in \GL(n,\C)$. We define
the compact quantum group $\cG = A_u(F)$ as follows.
\begin{itemize}
\item $C(\cG)$ is the universal \cst-algebra with generators
  $(U_{ij})$ and relations making $U = (U_{ij})$ and $F \overline{U}
  F^{-1}$ unitary elements of $\M_n(\C) \ot C(\cG)$, where $(\overline{U})_{ij} = U_{ij}^*$.
\item $\de(U_{ij}) = \sum_k U_{ik} \ot U_{kj}$.
\end{itemize}
\end{definition}

\begin{definition}
Let $F \in \GL(n,\C)$ satisfying $F \overline{F} = \pm 1$. We define
the compact quantum group $\cG = A_o(F)$ as follows.
\begin{itemize}
\item $C(\cG)$ is the universal \cst-algebra with generators
  $(U_{ij})$ and relations making $U = (U_{ij})$ a unitary element of
  $\M_n(\C) \ot C(\cG)$ and $U = F \overline{U} F^{-1}$, where $(\overline{U})_{ij} = U_{ij}^*$.
\item $\de(U_{ij}) = \sum_k U_{ik} \ot U_{kj}$.
\end{itemize}
\end{definition}

In both examples, the unitary matrix $U$ is a representation, called the \emph{fundamental representation}. The definition of $\cG=A_o(F)$ makes
sense without the requirement $F \overline{F} = \pm 1$, but the fundamental representation is irreducible if and only if $F \overline{F} \in \R 1$.

\begin{remark}
It is easy to classify the quantum groups $A_o(F)$. For $F_1,F_2 \in
\GL(2,\C)$ with $F_i \overline{F}_i = \pm 1$, we write $F_1 \sim F_2$
if there exists a unitary matrix $v$ such
that $F_1 = v F_2 v^t$, where $v^t$ is the transpose of $v$. Then,
$A_o(F_1) \cong A_o(F_2)$ if and only if $F_1 \sim F_2$.
It follows
that the $A_o(F)$ are classified up to isomorphism by $n$, the sign $F
\overline{F}$ and the eigenvalue list of $F^*F$ (see e.g.\ Section 5
of \cite{BDRV} where an explicit fundamental domain for the relation
$\sim$ is described).

If $F \in \GL(2,\C)$, we get up to isomorphism, the matrices
$$F^\pm_q = \begin{pmatrix} 0 & \sqrt{q} \\ \mp \frac{1}{\sqrt{q}} & 0
\end{pmatrix}$$
for $0 < q \leq 1$ with corresponding quantum groups $A_o(F^\pm_q) \cong \SU_{\pm q}(2)$.

{\it For the rest of the paper, we shall assume that $F \not\sim
  F^\pm_1$, which means that we do not deal with the classical group
  $\SU(2)$, neither with $\SU_{-1}(2)$.}
Generally speaking, our interest lies in $A_o(F)$ with $\dim F \geq 3$.
\end{remark}

The quantum groups $A_u(F)$ and $A_o(F)$ have been studied extensively
by Banica \cite{banica1,banica2}. In particular, he gave a complete
description of their representation theory. In the rest of the paper
we focus on $A_o(F)$. The following result is proven in
\cite{banica1}: it tells us that $A_o(F)$ has the same fusion rules as
the classical compact group $\SU(2)$. Observe however that the
dimension of the fundamental representation $U$ is $n$. Conversely, it
is easy to see that any compact quantum group with the same fusion
algebra as $\SU(2)$ is isomorphic to an $A_o(F)$.

\begin{theorem}
Let $F \in \GL(n,\C)$ and $F \overline{F} = \pm 1$. Let $\cG = A_o(F)$. One can identify $\Irred(\cG)$ with $\N$ in such a way that
$$x \ot y \cong |x-y| \oplus (|x-y|+2) \oplus \cdots \oplus (x+y) \;
,$$
for all $x,y \in \N$.
\end{theorem}

It is easy to check that $\dimq(1) = \Tr(F^*F)$. Take $0 < q < 1$ such that $\Tr(F^*F) = q + \frac{1}{q}$. Then,
$$\dimq(n) = [n+1]_q \quad\text{where}\quad [n]_q = \frac{q^n
  - q^{-n}}{q - q^{-1}} \; .$$
When there is no confusion, we do not write the index $q$ in the $q$-number $[n]_q$.

\begin{notation} \label{not.intertwiners}
Let $\cG = A_o(F)$ and $x,y \in \Irred(\cG)$. Whenever $z \in x \ot
y$, we denote by $V(x \ot y,z)$ an \emph{isometric} element in $\Mor(x
\ot y,z)$. Note that $V(x \ot y,z)$ is defined up to a number of
modulus $1$.

Whenever $z \in x \ot
y$, we denote by $p^{x \ot y}_z$ the unique orthogonal projection in
$\End(x \ot y)$ projecting onto the irreducible representation
equivalent with $z$.
\end{notation}

Throughout the paper, the letters $n,x,y,z,r,s$ are reserved to denote irreducible representations of $A_o(F)$ (i.e.\ natural numbers). The letters
$a,b,c,\ldots$ are used to denote elements of \cst-algebras. The capital letters $A$ and $B$ denote matrices.

\section{Solidity and the Akemann-Ostrand property} \label{sec.solid}

In \cite{ozawa}, Ozawa introduced the following remarkable definition. Recall that a von Neumann algebra is said to be diffuse if it does not contain
minimal projections.

\begin{definition}[N.\ Ozawa]
A von Neumann algebra $M$ is said to be \emph{solid} if $M \cap A'$ is injective for any \emph{diffuse} subalgebra $A \subset M$.
\end{definition}

A solid von Neumann algebra is necessarily finite. The following definition is a straightforward adaptation of solidity to arbitrary von Neumann
algebras and has been observed independently by D.\ Shlyakhtenko \cite{dima-personal}.

{\it From now on, we assume that von Neumann algebras have separable predual.}

\begin{definition}
A von Neumann algebra $M$ is said to be \emph{generalized solid} if $M \cap A'$ is injective for any \emph{diffuse} subalgebra $A \subset M$ for
which there exists a \emph{faithful normal conditional expectation} $E : M \recht A$.
\end{definition}

Several results in \cite{ozawa} can now be easily generalized. For the convenience of the reader, we give an overview of what we need in this paper.
The first result is immediate.

\begin{proposition} \label{prop.solidprime}
\begin{itemize}
\item A finite von Neumann algebra is generalized solid if and only if it is solid.
\item A subalgebra $M_1 \subset M$ of a generalized solid von Neumann algebra that admits a faithful normal conditional expectation $M \recht M_1$,
is again generalized solid.
\item A non-injective generalized solid factor $M$ is prime: if $M \cong M_1 \ot M_2$, then either $M_1$ or $M_2$ is a
type I factor.
\end{itemize}
\end{proposition}

The main result of \cite{ozawa} consists in deducing solidity from the \emph{Akemann-Ostrand property}. Recall the following definition from
\cite{ozawa}.

\begin{definition}
A von Neumann algebra $M \subset \B(H)$ is said to satisfy the \emph{Akemann-Ostrand property}, if there exist unital weakly dense \cst-subalgebras
$B \subset M$, $C \subset M'$ such that $B$ is locally reflexive and the $^*$-homomorphism
$$B \otalg C \recht \frac{\B(H)}{\cK(H)} : \sum_{i=1}^n b_i \ot c_i \mapsto \sum_{i=1}^n \pi(b_i c_i)$$
extends continuously to $B \otmin C$. Here $\pi$ denotes the quotient map $\B(H) \recht \B(H)/\cK(H)$.
\end{definition}

Theorem 6 in \cite{ozawa} has the following generalization.

\begin{theorem} \label{thm.aosolid}
A von Neumann algebra $M \subset \B(H)$ satisfying the Akemann-Ostrand property is generalized solid.
\end{theorem}
\begin{proof}
One follows almost line by line the proof in \cite{ozawa}, only paying attention that there are conditional expectations everywhere since they do not
exist automatically on von Neumann subalgebras (contrary to the finite case). Suppose $A \subset M$ is diffuse and $E : M \recht A$ a faithful normal
conditional expectation. Choose on $A$ a faithful state $\vphi$ such that the centralizer algebra $A^\vphi$ has a diffuse abelian subalgebra $A_0
\subset A^\vphi$. This is indeed possible, the most difficult case of $A$ a type III$_1$ factor being dealt with in \cite{CS}, Corollary 8.

Write $\psi = \vphi E$. Since there is a $\psi$-preserving conditional expectation $M \cap A_0' \recht M \cap A'$, it is sufficient to show that $M
\cap A_0'$ is injective. Because there is a unique $\psi$-preserving conditional expectation $M \recht M \cap A_0'$, the proof of \cite{ozawa}
applies literally.
\end{proof}

On the level of compact quantum group \cst-algebras, we have the
following version of the Akemann-Ostrand property.

\begin{definition}
Let $\cG$ be a compact quantum group. We say that $\cG$ satisfies the
\emph{Akemann-Ostrand property} if the $^*$-homomorphism
$$C(\cG)\red \otalg C(\cG)\red \recht
\frac{\B(L^2(\cG))}{\cK(L^2(\cG))} : \sum_{i=1}^n a_i \ot b_i \mapsto
\sum_{i=1}^n \pi(\lambda(a_i) \rho(b_i))$$
extends continuously to $C(\cG)\red \otmin C(\cG)\red$. Here $\pi$
denotes the quotient map $\B(L^2(\cG)) \recht \frac{\B(L^2(\cG))}{\cK(L^2(\cG))}$.
\end{definition}

Obviously, if $\cG$ satisfies the Akemann-Ostrand property and if $C(\cG)\red$ is locally reflexive, the von Neumann algebra $C(\cG)\vnalg$ satisfies
the Akemann-Ostrand property as well and is, by Theorem \ref{thm.aosolid}, a generalized solid von Neumann algebra.

\section{Boundary and boundary action for the dual of $A_o(F)$} \label{sec.boundary}

Fix $F \in \GL(n,\C)$ with $F \overline{F} = \pm 1$. Put $\cG =
A_o(F)$. Recall that we assume that $\cG \not\cong SU_{\pm
  1}(2)$. Recall that we identify $\Irred(\cG) = \N$ and that we use
the letters $n,x,y,z,r,s$ to denote irreducible representations of $\cG$.

We shall introduce a boundary for $\cGh$, inspired by the construction
of the boundary of a free group by adding infinite reduced words. So,
we first define a \emph{compactification} of $\cGh$, which will be a
unital \cst-algebra $\cB$ such that
$$c_0(\cGh) \subset \cB \subset \ell^\infty(\cGh) \; .$$
The boundary $\cB_\infty$ is then defined as $\cB_\infty = \cB / c_0(\cGh)$.

We show that the comultiplication $\deh$ yields by restriction and passage to the quotient $\cB_\infty = \cB/c_0(\cGh)$
\begin{itemize}
\item an action of $\cGh$ on $\cB_\infty$ on the left;
\item the trivial action of $\cGh$ on $\cB_\infty$ on the right.
\end{itemize}
In the next section we shall introduce the notion of an amenable action and prove that the boundary action is amenable.

If one compactifies a free group $\Gamma$ by adding infinite words, a continuous function on this compactification is an element of
$\ell^\infty(\Gamma)$ whose value in a long word of $\Gamma$ essentially only depends on the beginning part of that word. In order to give somehow
the same kind of definition for $\cGh$ the dual of $A_o(F)$, we need to compare the values that an element of $\ell^\infty(\cGh)$ takes in two
different irreducible representations. So, we should compare matrices in $\B(H_x)$ and $\B(H_y)$ for $x,y \in \Irred(\cGh)$. To do so, we use the
following linear maps.

\begin{definition} \label{def.psi}
Let $x,y \in \N$. We define unital completely positive maps
$$\psi_{x+y,x} : \B(H_x) \recht \B(H_{x+y}) : \psi_{x+y,x}(A) = V(x \ot y,x+y)^* (A \ot 1) V(x \ot y,x+y) \; .$$
Recall that we have chosen isometric intertwiners $V(x \ot y,z) \in \Mor(x \ot y,z)$.
\end{definition}

\begin{proposition}
The maps $\psi_{x+y,x}$ form an inductive system of completely positive maps.
\end{proposition}
\begin{proof}
Since $\Mor(x+y+z,x \ot y \ot z)$ is one-dimensional, we have that
$$
(V(x \ot y,x+y) \ot 1)V((x+y) \ot z,x+y+z) = (1 \ot V(y \ot z,y+z))V(x \ot (y+z),x+y+z) \mod \T \; .
$$
So, we are done.
\end{proof}

\begin{notation}
We define
$$\psi_{\infty,x} : \B(H_x) \recht \ell^\infty(\cGh) : \psi_{\infty,x}(A)p_y = \begin{cases} \psi_{y,x}(A) \;\;\text{if}\;\; y \geq x \\ 0 \;\;\text{else}\end{cases} \; .$$
We use the same notation for the map $\psi_{\infty,x} : \B(H_x) \recht \frac{\ell^\infty(\cGh)}{c_0(\cGh)}$. Recall that $p_x$ denotes the minimal
central projection in $\ell^\infty(\cGh)$ associated with $x \in \Irred(\cG)$.
\end{notation}

\begin{proposition} \label{prop.compactif}
Define
$$\cB_0 = \{ a \in \ell^\infty(\cGh) \mid \;\;\text{there exists}\;\; x \;\;\text{such that}\;\; ap_y = \psi_{y,x}(ap_x) \;\;\text{for all}\;\; y \geq x \} \;
.$$ The norm closure of $\cB_0$ is a unital \cst-subalgebra of $\ell^\infty(\cGh)$ containing $c_0(\cGh)$. We denote it by $\cB$. The \cst-algebra
$\cB$ is nuclear.
\end{proposition}

\begin{proof}
It follows from \eqref{eq.best1} that there exists a constant $C$ such that
$$\| [(\psi_{x+y,x}(A) \ot 1), p^{(x+y)\ot z}_{x+y+z}] \| \leq C q^{y} \|A\|$$
for all $x,y,z$ and $A \in \B(H_x)$. Hence,
$$\| \psi_{x+y+z,x+y}\bigl( \psi_{x+y,x}(A) B \bigr) - \psi_{x+y+z,x}(A) \psi_{x+y+z,x+y}(B) \| \leq C q^{y} \|A\| \, \|B\|$$
for all $x,y,z$, $A \in \B(H_x)$ and $B \in \B(H_{x+y})$. We easily conclude that the norm closure $\cB$ is a unital \cst-subalgebra of
$\ell^\infty(\cGh)$. It obviously contains $c_0(\cGh)$.

Define $B_n = \bigoplus_{x=0}^n \B(H_x)$ and define $\mu_n : \cB \recht B_n$ by restriction and $\gamma_n : B_n \recht \cB$ by the formula
$\gamma_n(a)p_x = ap_x$ if $x \leq n$ and $\gamma_n(a)p_x = \psi_{x,n}(ap_n)$ if $x \geq n$. Then $\gamma_n (\mu_n(a)) \recht a$ for all $a \in \cB$
and the nuclearity of $\cB$ is proven.
\end{proof}

\begin{notation}
The comultiplication $\deh$ yields a (left) action of $\cGh$ on the \cst-algebra $\ell^\infty(\cGh)$ that we denote by
$$\be : \ell^\infty(\cGh) \recht \M(c_0(\cGh) \ot \ell^\infty(\cGh)) \; .$$
\end{notation}

\begin{proposition}
We have $\be(\cB) \subset \M(c_0(\cGh) \ot \cB)$ and as such $\be$ is an action of $\cGh$ on $\cB$.
\end{proposition}

\begin{proof}
It suffices to show that $(p_x \ot 1) \be(a) \in \B(H_x) \ot \cB$ for all $a \in \cB$ and $x \in \N$. Take $a=\psi_{\infty,r}(A)$. Take $y \geq x+r$
and take $z$. Then,
$$(p_x \ot p_{y+z}) \be(x) = \sum_{s \in x \ot y} V(x \ot (y+z),s+z) \, (a p_{s+z}) \, V(x \ot (y+z),s+z)^* \; .$$
Fix an $s \in x \ot y$. Observe that $ap_{s+z} = \psi_{s+z,s}(\psi_{s,r}(A))$. Using \eqref{eq.best3}, we get that
\begin{multline*}
\| V(x \ot (y+z),s+z) \, (a p_{s+z}) \, V(x \ot (y+z),s+z)^* - (\id \ot \psi_{y+z,y})\bigl( V(x \ot y,s) \psi_{s,r}(A) V(x \ot y,s)^* \bigr) \|
\\ \leq C q^{-x+y} \|A\| \; .
\end{multline*}
We keep $x$ and $A \in \B(H_r)$ fixed. Choose $\eps > 0$. Take $y$ such that $(x+1)C q^{-x+y} \|A\| < \eps$. Since there are less then $x+1$
irreducible components of $x \ot y$, the computation above shows that $(p_x \ot 1) \be(a)$ is at distance at most $\eps$ of
$$(\id \ot \psi_{\infty,y})\Bigl( \sum_{s \in x \ot y} V(x \ot y,s) \psi_{s,r}(A) V(x \ot y,s)^* \Bigr) \; .$$
Hence, $(p_x \ot 1) \be(a) \in \B(H_x) \ot \cB$.
\end{proof}

\begin{definition}
We define $\cB_\infty := \cB/c_0(\cGh)$ and we still denote by $\be$ the action of $\cGh$ on $\cB_\infty$.
\end{definition}

As it is the case for the action of a free group on its boundary, we prove that the action by right translations on $c_0(\cGh)$ extends to an action
on $\cB$ that becomes the trivial action on $\cB_\infty$. The precise statement is as follows.

\begin{proposition} \label{prop.higson}
Consider the (right) action $\gamma : \ell^\infty(\cGh) \recht \M(\ell^\infty(\cGh) \ot c_0(\cGh))$ of $\cGh$ on the \cst-algebra $\ell^\infty(\cGh)$
by right translations. For all $a \in \cB$ and all $x$, we have
$$\bigl(\gamma(a) - a \ot 1 \bigr)(1 \ot p_x) \in c_0(\cGh) \ot \B(H_x) \;
.$$ Hence, $\gamma$ becomes the trivial action on $\cB_\infty$.
\end{proposition}

\begin{proof}
Suppose $a = \psi_{\infty,x}(A)$ for $A \in \B(H_x)$. Fix a $z$ and take $y \geq z$. Using \eqref{eq.best-a3}, we get a constant $C$ such that,
\begin{align*}
\deh & (\psi_{\infty,x}(A))  (p_{x+y} \ot p_z) = \sum_{s \in y \ot z} V((x+y) \ot z,x+s)\psi_{x+s,x}(A) V((x+y) \ot z,x+s)^* \\ & \approx \sum_{s \in
y \ot z} \bigl(V(x
\ot y,x+y)^* \ot 1 \bigr) (A \ot p^{y \ot z}_{s}) \bigl(V(x \ot y,x+y) \ot 1 \bigr) \quad\text{with error}\;\; \leq (z+1) C q^{-z+y} \\
&= \psi_{x+y,x}(A) \ot p_z \; .
\end{align*}
If we keep fixed $A$ and let $y \recht \infty$, the conclusion follows.
\end{proof}

For later use, we prove the following lemma. The only interest at this
point, is that it shows that $\cB_\infty$ is non-trivial: the maps
$\psi_{\infty,x} : \B(H_x) \recht \cB_\infty$ are injective.

\begin{lemma} \label{lem.injective-psi}
There exists a constant $D > 0$, only depending on $q$, such that
$$D \|A\|_{\psi_x} \leq \|\psi_{x+y,x}(A)\|_{\psi_{x+y}} \leq \|A\|_{\psi_x}$$
for all $x,y$ and $A \in \B(H_x)$.
\end{lemma}
\begin{proof}
Consider $\B(H_x)$ as a Hilbert space using the state $\psi_x$. Then,
$$B(H_x) \recht H_x \ot H_x : A \mapsto (A \ot 1)t_x$$
is a unitary operator. Using the notation $D(x,y) = [x+1] [y+1] [x+y+1]^{-1}$ and \eqref{eq.relative}, we know that $t_{x+y}$ equals, up to a number
of modulus one,
$$D(x,y)^{1/2} \bigl(V(x \ot y,x+y)^* \ot V(y \ot x,x+y)^*\bigr) (1 \ot t_y \ot 1)t_x \; .$$
Using Lemma \ref{lem.encore-une}, we get a constant $D > 0$ such that
\begin{align*}
\|\psi_{x+y,x}(A)\|_{\psi_{x+y}} &= \| (\psi_{x+y,x}(A) \ot 1)t_{x+y} \| \\
&= D(x,y)^{1/2} \|(V(x \ot y,x+y)^* \ot 1)(A \ot 1 \ot 1)(p^{x \ot y}_{x+y} \ot V(y \ot x,x+y)^*)(1 \ot t_y \ot 1)t_x \| \\
&= D(x,y)^{1/2} \|(V(x \ot y,x+y)^* \ot 1)(A \ot 1 \ot 1)(1 \ot 1 \ot V(y \ot x,x+y)^*)(1 \ot t_y \ot 1)t_x \| \\
&= D(x,y)^{1/2} \|(V(x \ot y,x+y)^* \ot V(y \ot x,x+y)^*)(1 \ot t_y \ot 1) \; (A \ot 1)t_x \| \\
&\geq D \|(A \ot 1)t_x \| = D \|A\|_{\psi_x} \; .
\end{align*}
So, we are done.
\end{proof}

\section{Amenability of the boundary action and the Akemann-Ostrand
  property} \label{sec.amenable}

We introduce the notion of an \emph{amenable action} of a discrete
quantum group on a unital \cst-algebra. We prove that for $\cG = A_o(F)$, the
action of $\cGh$ on its boundary $\cB_\infty$ as introduced in the
previous section, is amenable. We then deduce the exactness of
$C(\cG)\red$ and the Akemann-Ostrand property.

In the following definition, we make use of the representation
$$\rhoh \lambdah : \ell^\infty(\cGh) \ovt \ell^\infty(\cGh) \recht \B(L^2(\cG))
: (\rhoh \lambdah)(a \ot b) \xi_x = (bp_x \ot ap_{\ox})\xi_x \quad\text{for
  all}\quad \xi_x \in H_x \ot H_{\ox} \; .$$
Recall that $\cVh$ denotes the right regular representation of $\cGh$.

\begin{definition} \label{def.amenableaction}
Let $\be : \cB \recht \M(c_0(\cGh) \ot \cB)$ be an action of the discrete quantum group $\cGh$ on a unital \cst-algebra $\cB$. We say that $\be$ is
\emph{amenable} if there exists a sequence $\xi_n \in L^2(\cG) \ot \cB$ satisfying
\begin{itemize}
\item $\xi_n^* \xi_n \recht 1$ in $\cB$,
\item for all $x$, $\bigl\| \bigl( (\id \ot \be)(\xi_n) - \cVh_{12} (\xi_n)_{13} \bigr) (1 \ot p_x \ot 1) \bigr\| \recht 0$.
\item $(\rhoh\lambdah\deh \ot \id)\be(a) \xi_n = \xi_n a$ for all $n$ and all $a \in \cB$,
\end{itemize}
\end{definition}

\begin{remark}
It is clear that $\cGh$ is amenable if and only if the trivial action on $\C$ is amenable if and only if every action is amenable.

The first two conditions in the previous definition are natural: the $\xi_n$ are approximately equivariant unit vectors. The third condition might
seem mysterious. But, already the definition of an amenable action of a discrete group on a unital \cst-algebra involves an extra condition: see
Th{\'e}or{\`e}me 3.3 in \cite{AD} where positive definite functions take values in the center. In the quantum setting, this centrality condition is replaced
by the third condition in the definition above, and reduces to centrality in the case where $\cGh$ is a discrete group. Indeed, in that case
$(\rhoh\lambdah)\deh$ is the co-unit $\epsh$ and the condition above reads $(1 \ot a)\xi_n = \xi_n a$ for all $a \in \cB$ and all $n$.
\end{remark}

\begin{notation}
In this section, we write $H$ for the Hilbert space $L^2(\cG)$.
\end{notation}

\begin{proposition} \label{prop.fullreduced}
Let $\be : \cB \recht \M(c_0(\cGh) \ot \cB)$ be an amenable action of a discrete quantum group $\cGh$ on a unital \cst-algebra $\cB$. Then, the
natural homomorphism $\cGh \ltimesfull \cB \recht \cGh \ltimesred \cB$ is an isomorphism. If moreover $\cB$ is nuclear, $\cGh \ltimes \cB$ is
nuclear, the reduced \cst-algebra $C(\cG)\red$ is an exact \cst-algebra and $\cGh$ is an exact quantum group.
\end{proposition}

\begin{proof}
Let $(\theta,X)$ be a covariant representation of $\be$ on the Hilbert
space $K$. Define bounded linear maps
$$v_n : K \recht H \ot K : v_n \eta = (\lambdah \ot \id)(X)(\id \ot \theta)(\xi_n) \eta \; .$$
We shall prove that the $v_n$ approximately intertwine the covariant
representation $(\theta,X)$ with a regular covariant
representation. First observe that, for all $a \in \cB$,
$$v_n \theta(a) = (\lambdah \ot \id)(X) \bigl( (\rhoh\lambdah \ot \theta)(\id
\ot \be)\be(a) \bigr) (\id \ot \theta)(\xi_n) \; .$$
For all $a \in c_0(\cGh)$ and $b \in \cB$, we have
\begin{align*}
(\lambdah \ot \id)(X) & \bigl((\rhoh\lambdah \ot \theta)(\id \ot \beta)(a \ot
b)\bigr) = (\lambdah \ot \id)(X) (\rhoh(a) \ot 1) (\lambdah \ot \id)\beta(b)
\\
&= (\rhoh(a) \ot 1) (1 \ot \theta(b))(\lambdah \ot \id)(X) = (\rhoh \ot
\theta)(a \ot b)(\lambdah \ot \id)(X) \; .
\end{align*}
Hence,
$$v_n \theta(a) = (\rhoh \ot \theta)\be(a) v_n$$
for all $a \in \cB$.

Next, observe that
$$((1 \ot v_n)X (p_x \ot 1))_{213} = (\lambdah \ot \id)(X)_{13} X_{23}
(\id \ot \id \ot \theta)\bigl((\id \ot \be)(\xi_n) (1 \ot p_x \ot 1)\bigr) \; ,$$ while
$$( \cVh_{21} (1 \ot v_n) (p_x \ot 1))_{213} = (\lambdah \ot \id)(X)_{13}
X_{23} (\id \ot \id \ot \theta)\bigl( \cVh_{12} (\xi_n)_{13} (1 \ot p_x \ot 1)\bigr) \; .$$
The condition in \ref{def.amenableaction} yields that
$$(p_x \ot 1 \ot 1)\bigl((1 \ot v_n)X - \cVh_{21} (1 \ot v_n)\bigr)
\recht 0$$
for all $x \in \Irred(\cG)$. So, we have shown that the $v_n$
approximately intertwine $(\theta,X)$ with the regular covariant
representation $((\rhoh \ot \theta)\be, \cVh_{21})$.
Since
$v_n^* v_n \recht 1$ in the norm topology, this shows that $\cGh \ltimesfull \cB \recht \cGh \ltimesred \cB$ is an isomorphism.

In order to show that $\cGh \ltimes \cB$ is nuclear when $\cB$ is nuclear, it suffices to observe that the action $\be \ot \id$ of $\cGh$ on $\cB \ot
D$ is amenable when $\be$ is amenable. As a subalgebra of a nuclear
\cst-algebra, the reduced \cst-algebra $C(\cG)\red$ is exact. The
exactness of $\cGh$ follows from Proposition \ref{prop.quantumexact}.
\end{proof}

We now fix $F \in \GL(n,\C)$ with $F \overline{F} = \pm 1$ and take for the rest of this section $\cG = A_o(F)$. We still have our standing
assumption that $\cG \not\cong SU_{\pm 1}(2)$.

\begin{theorem}
Let $\cG=A_o(F)$. The boundary action of $\cGh$ on $\cB_\infty$, constructed in the previous section, is amenable.
\end{theorem}

\begin{proof}
Consider the unit vector $\mu := \cVh (\xi_0 \ot 1) \in \cL(c_0(\cGh),H \ot c_0(\cGh))$ as well as the unit vectors
$$\mu_x = \mu p_x \in H \ot \B(H_x) \; .$$
Note that $\mu = (\Lambdahr \ot \id)\deh(p_\eps)$.
Observe that $\mu_x \in \lambdah(p_{\ox})H \ot \B(H_x)$, which implies that the vectors $\mu_x$ are mutually orthogonal.

Define $\xi_n \in H \ot \cB_\infty$ by the formula $\displaystyle \xi_n = \frac{1}{\sqrt{n+1}} \sum_{y=0}^n (\id \ot \psi_{\infty,y})(\mu_y)$. We
claim that $\xi_n^* \xi_n \recht (1-q^{2}) 1$. Consider, for a fixed $x$, the vector $\eta_x \in H \ot \cB_\infty$ given by $\eta_x = (\id \ot
\psi_{\infty,x})(\mu_x)$. Since $\mu_x = (\Lambdahr \ot \id)(p^{\ox \ot x}_\eps)$, we get
\begin{align*}
(\mu_x^* \ot 1)(1 \ot p^{x \ot y}_{x+y})(\mu_x \ot 1) &= [x+1]^2
(\vphi_{\ox} \ot \id \ot \id)\bigl( (p^{\ox \ot x}_\eps \ot 1)(1 \ot
p^{x \ot y}_{x+y})(p^{\ox \ot x}_\eps \ot 1)\bigr) \\ &= [x+1]^2
\bigl( (\vphi_{\ox} \ot \id)(p^{\ox \ot x}_\eps) \ot 1\bigr) \bigl(1
\ot (\vphi_x \ot \id)(p^{x \ot y}_{x+y})\bigr) \\ &=
\frac{[x+y+1]}{[x+1] [y+1]} 1 \ot 1 \; .
\end{align*}
It follows that $\eta_x^* \eta_x = \frac{q^{-x}}{[x+1]} 1$ in $\cB_\infty$. Since $\frac{q^{-x}}{[x+1]} \recht 1 - q^{2}$ when $x \recht \infty$, the
claim is proven.

In order to verify that $(\rhoh\lambdah\deh \ot \id)\be(a) \xi_n = \xi_n a$
for all $n$ and all $a \in \cB_\infty$, it is sufficient to check that
\begin{equation}\label{eq.againsufficient}
(\id \ot \psi_{x+y,x})(\mu_x) (ap_{x+y}) = \bigl((1 \ot
p_{x+y})(\rhoh\lambdah\deh \ot \id)\be(a) \bigr) (\id \ot
\psi_{x+y,x})(\mu_x)
\end{equation}
for all $a \in \ell^\infty(\cGh)$ and all $x,y \in
\Irred(\cG)$. Observe that the right hand side of
\eqref{eq.againsufficient} equals
\begin{equation} \label{eq.weereentje}
(1 \ot V(x \ot y,x+y)^*) \bigl( (1 \ot p_x \ot p_y) (\rhoh\lambdah \ot \id \ot \id)\deh^{(3)}(a) (\mu \ot 1) \bigr) V(x \ot y,x+y) \; .
\end{equation}
But, for all $a,b \in c_0(\cGh)$, we have
\begin{align*}
(\rhoh\lambdah \ot \id)(a \ot \deh(b)) \mu &= (\rhoh(a) \ot 1) \, (\lambdah \ot
\id)\deh(b) \, \cVh (\xi_0 \ot 1) = (\rhoh(a) \ot 1) \cVh (\lambdah(b)\xi_0
\ot 1) \\ &= \cVh (\rhoh \ot \id)\dehop(a) (\xi_0 \ot 1) \epsh(b) = \mu
a \epsh(b) \; .
\end{align*}
It follows that \eqref{eq.weereentje} equals
$$(1 \ot V(x \ot y,x+y)^*) (\mu_x \ot 1) \bigl( (p_x \ot p_y) \deh(a)
\bigr) V(x \ot y,x+y) = (\id \ot \psi_{x+y,x})(\mu_x) (ap_{x+y}) \;
.$$
This proves \eqref{eq.againsufficient}.
We then come to the crucial approximate equivariance condition for
$\xi_n$. First observe that $\cVh_{12} \mu_{13} = \cVh_{12} \cVh_{13}
(\xi_0 \ot 1 \ot 1) = (\id \ot \deh)(\mu)$. Hence,
for all $x$ and all $y \geq n$,
\begin{align*}
\cVh_{12} (\xi_n)_{13} (p_x \ot p_y) &= \frac{1}{\sqrt{n+1}} \sum_{s=0}^n (\id \ot \id \ot \psi_{y,s})\bigl( \cVh_{12} \mu_{13} (p_x \ot p_s) \bigr) \\
&=\frac{1}{\sqrt{n+1}} \sum_{s=0}^n (\id \ot \id \ot \psi_{y,s})\bigl( (\id \ot \deh)(\mu) (p_x \ot p_s) \bigr) \; .
\end{align*}
Take $n \geq K \geq x$ and $y \geq n+x$. We now write equalities up to an error term, that we estimate using the norm of the \cst-module $H \ot
\B(H_x) \ot \B(H_y)$.
\begin{align*}
\cVh_{12} (\xi_n)_{13} (p_x \ot p_y) & \approx \frac{1}{\sqrt{n+1}} \sum_{s=K}^n (\id \ot \id \ot \psi_{y,s})\bigl( (\id \ot \deh)(\mu) (p_x \ot p_s)
\bigr) \quad\text{with error}\;\; \leq
\frac{K}{\sqrt{n+1}} \\
& = \frac{1}{\sqrt{n+1}} \sum_{s=K}^n \sum_{z \in x \ot s} (\id \ot \id \ot \psi_{y,s})\bigl( (1 \ot V(x \ot s,z)) \mu_z V(x \ot s,z)^* \bigr) =
(\star) \; .
\end{align*}
It follows from \eqref{eq.best3} that there exists a constant $C$ such that
\begin{align*}
\bigl\|  (\id & \ot \id  \ot \psi_{y,s})\bigl( (1 \ot V(x \ot s,z)) \mu_z V(x \ot s,z)^* \bigr) \\ & \hspace{2cm} - \bigl( 1 \ot V(x \ot y,z+y-s)
\bigr) (\id \ot \psi_{z+y-s,z})(\mu_z) V(x \ot
y,z+y-s)^* \bigr\| \\
& \leq 2 C q^{-x+s} \leq 2 C q^{-x+K} \; .
\end{align*}
Observe now that in the sum $(\star)$, for a given $z$ there are less then $x+1$ corepresentations $s$ such that $z \in x \ot s$. Moreover, in the
sum $(\star)$, $z$ ranges from $K-x$ to $n+x$ and we have $\mu_z \in \lambdah(p_z)H \ot \B(H_z)$, with the $\lambdah(p_z)H$ orthogonal for different
$z$. We also observe that $z \in x \ot s$ if and only if $z+y-s \in x \ot y$ and conclude that
\begin{align*}
(\star) &\approx \frac{1}{\sqrt{n+1}} \sum_{s=K}^n \sum_{z \in x \ot
  s} \bigl( 1 \ot V(x \ot y,z+y-s)\bigr) (\id \ot
  \psi_{z+y-s,z})(\mu_z) V(x \ot y,z+y-s)^* \\
&\hspace{11cm}\text{with error}\;\; \leq 4(x+1) C q^{-x+K} \\
& = \frac{1}{\sqrt{n+1}} \sum_{s=K}^n \sum_{r \in  x \ot y} ( 1 \ot V(x \ot y,r)) (\id \ot \psi_{r,r-y+s})(\mu_{r-y+s}) V(x \ot y,r)^* \\
& = (\id \ot \deh)(\eta)(p_x \ot p_y) \; ,
\end{align*}
where $\eta \in H \ot c_0(\cGh)$ is given by $\displaystyle \eta p_r = \frac{1}{\sqrt{n+1}} \sum_{s=K}^n (\id \ot \psi_{r,r-y+s})(\mu_{r-y+s})$
whenever $y-x \leq r \leq y+x$ and $\eta p_r = 0$ elsewhere. When $y-x \leq r \leq y+x$,
$$\bigl\| \eta p_r - \xi_n p_r \bigr\| \leq \frac{2x+K}{\sqrt{n+1}} \; .$$
Since the expression $(\id \ot \deh)(\eta)(p_x \ot p_y)$ only takes into account the values of $\eta p_r$ for $y-x \leq r \leq y+x$, it follows that
$$(\id \ot \deh)(\eta)(p_x \ot p_y) \approx (\id \ot \deh)(\xi_n)(p_x \ot p_y) \quad\text{with error}\;\; \leq \frac{2x+K}{\sqrt{n+1}} \; .$$
We finally conclude that
$$\bigl\| \bigl(\cVh_{12} (\xi_n)_{13} - (\id \ot \be)(\xi_n) \bigr) (1 \ot p_x \ot 1) \bigr\| \leq \frac{2x+2K}{\sqrt{n+1}} + 4 C (x+1)q^{-x+K} \; .$$
Given $x$, we first take $K$ such that $4 C (x+1) q^{-x+K}$ is small. We then take $n$ such that $\frac{2x+2K}{\sqrt{n+1}}$ is small. As such, we
have shown the amenability of the action $\be$. Indeed, it suffices to replace $\xi_n$ by $\bigl(\frac{1}{1 - q^2}\bigr)^{1/2} \xi_n$.
\end{proof}

\begin{remark}
The same proof shows that the action of $\cGh$ on $\ell^\infty(\cGh)/c_0(\cGh)$ by left translation is an amenable action. But, since
$\ell^\infty(\cGh)/c_0(\cGh)$ is non-nuclear (even non-exact), we really need the amenability of the action on the nuclear \cst-algebra $\cB_\infty$
to show, e.g., the exactness of $C(\cG)\red$.
\end{remark}

We deduce exactness and the Akemann-Ostrand property from the
amenability of the boundary action. Note that an independent proof of
the Akemann-Ostrand property has been given by the second author in \cite{Ver}.

\begin{corollary} \label{cor.exactAO}
Let $\cG = A_o(F)$. Then, $C(\cG)\red$ is exact and $\cG$ satisfies the Akemann-Ostrand property.
\end{corollary}

\begin{proof}
The exactness of $C(\cG)\red$ follows from Proposition
\ref{prop.fullreduced}. Put $H = L^2(\cG)$. Consider the left-right representation
$$\lambda \rho : C(\cG)\red \otmax C(\cG)\red \recht \frac{\B(H)}{\cK(H)} \; .$$
We have to show that this homomorphism factorises through $C(\cG)\red \otmin C(\cG)\red$. But, we also have the homomorphism $\cB_\infty \recht
\frac{\B(H)}{\cK(H)}$. It follows from Proposition \ref{prop.higson} that $\rho(C(\cG)\red)$ and $\cB_\infty$ commute in $\frac{\B(H)}{\cK(H)}$.
Hence, we get a homomorphism
$$(\cGh \ltimes \cB_\infty) \otmax C(\cG)\red \recht \frac{\B(H)}{\cK(H)} \; .$$
Since $\cGh \ltimes \cB_\infty$ is nuclear, the left hand side equals $(\cGh \ltimes \cB_\infty) \otmin C(\cG)\red$ and we are done.
\end{proof}

Combining with Theorem \ref{thm.aosolid}, we get the following result.

\begin{corollary} \label{cor.solid}
Let $\cG = A_o(F)$ and denote $M = C(\cG)\vnalg$. Then, $M$ is a generalized solid von Neumann algebra.
\end{corollary}

\section{Probabilistic interpretations of the boundary $\cB_\infty$} \label{sec.probabilistic}

A natural setting where boundaries of discrete groups appear is by
considering (invariant) random walks on the group. One associates to
such a random walk a Poisson boundary, which is a probability space,
and a Martin boundary, which comes from a bona fide compactification
of the group.

Both notions of Poisson boundary and Martin boundary have been
generalized to random walks on discrete quantum groups, see
\cite{B,C,I,INT,NT}.

In this section, we show that the Martin boundary for the dual of
$A_o(F)$ is naturally isomorphic with the boundary $\cB_\infty$
constructed above. Moreover, the Poisson boundary is isomorphic with the von
Neumann algebra generated by $\cB_\infty$ in the GNS construction of a
natural harmonic state on $\cB_\infty$. Bounded harmonic elements of
$\ell^\infty(\cGh)$ are written with a Poisson integral
formula. Note in this respect that a theorem establishing the
link between Martin and Poisson boundary for general discrete quantum
groups has not yet been established, see \cite{NT}.

\begin{notation}
Recall the states $\vphi_x$ and $\psi_x$ introduced in Notation \ref{not.states}. For every probability measure $\mu$ on $\Irred(\cG)$, we consider
the states
$$\psi_\mu = \sum_x \mu(x) \psi_x \quad\text{and}\quad \vphi_\mu = \sum_x \mu(x) \vphi_x \; .$$
Associated with these states, are the Markov operators
$$P_\mu = (\vphi_\mu \ot \id)\deh \quad\text{and}\quad Q_\mu = (\id \ot \psi_\mu)\deh \; .$$
\end{notation}

Note that a state $\om$ is of the form $\psi_\mu$ if and only if the Markov operator $(\id \ot \om)\deh$ preserves the center of $\ell^\infty(\cGh)$
(see e.g. Proposition 2.1 in \cite{NT}). Also note that we have a convolution product $\mu \star \nu$ on the measures on $\Irred(G)$, such that
$\psi_{\mu \star \nu} = \psi_\mu \star \psi_\nu$ and $\vphi_{\mu \star \nu} = \vphi_\mu \star \vphi_\nu$.

The operators $P_\mu$ and $Q_\mu$ are the Markov operators associated with a \emph{quantum random walk}. Their restriction to the center of
$\ell^\infty(\cGh)$ yields a Markov operator for a classical random
walk on the state space $\Irred(\cG)$, with transition probabilities
$p(x,y)$ and $n$-step transition probabilities $p_n(x,y)$
given by
\begin{equation}\label{eq.nstep}
p_x p(x,y) = p_x Q_\mu(p_y) \; , \quad p_x p_n(x,y) = p_x Q_\mu^n(p_y)
\; .
\end{equation}
Note that $p_n(e,y) = \mu^{\star n}(y) = \psi_\mu^{\star
n}(p_y)$.

\begin{definition}
The probability measure $\mu$ on $\Irred(\cG)$ is said to be transient if $\sum_{n=0}^\infty p_n(x,y) < \infty$ for all $x,y \in \Irred(\cG)$.
\end{definition}

Contrary to the case of random walks on discrete groups, probability measures on $\Irred(\cG)$ are very often transient, see Proposition 2.6 in
\cite{NT}. In particular, if $\cG=A_o(F)$ with $\cG \not\cong \SU(2),\SU_{-1}(2)$, every probability measure not concentrated in $0$ is transient.

\subsection*{Poisson boundary}

\begin{definition}
For any probability measure $\mu$ on $\Irred(\cG)$, we define $$H^\infty(\cGh,\mu) = \{a \in \ell^\infty(\cGh) \mid Q_\mu(a) = a \} \; .$$
\end{definition}
The weakly closed linear space $H^\infty(\cGh,\mu)$ is in fact a von Neumann algebra, with product given by
$$a \cdot b = \slim_{n \recht \infty} Q_\mu^n(ab) \; .$$
Remark that the Poisson boundary has a natural interpretation as a
relative commutant in the study of infinite tensor product actions
\cite{I,Vactions}.

\begin{terminology}
The support of a measure $\mu$ on $\Irred(\cG)$ is denoted by $\supp \mu$. We say that $\mu$ is generating, if
$$\Irred(\cG) = \bigcup_{n = 1}^\infty \supp(\mu^{\star n}) \; .$$
\end{terminology}

The restriction of $\epsh$ to $H^\infty(\cGh,\mu)$ defines a normal state on $H^\infty(\cGh,\mu)$. This state is faithful when $\mu$ is generating.

From now on, we fix $\cG = A_o(F)$ for a given matrix $F$ satisfying $F \overline{F} = \pm 1$.

Since the fusion rules of $A_o(F)$ are abelian, we know from Proposition 1.1 in \cite{INT} that $H^\infty(\cGh,\mu)$ does not depend on the choice of
a \emph{generating} measure $\mu$. Moreover, a measure $\mu$ on $\N = \Irred(\cG)$, is generating if and only if its support contains an odd number.

The aim of this section is to define a \emph{harmonic measure} (i.e.\ a state) on the boundary $\cB_\infty$ and to write every harmonic function
(i.e.\ element of $H^\infty(\cGh,\mu)$) as an integral with respect to the harmonic measure.

\begin{proposition}
The formula
$$\om(a) = \lim_{n \recht \infty} \psi_n(a)$$
yields a well defined state on $\cB$ and $\om(c_0(\cGh)) = \{0\}$. The resulting state on $\cB_\infty = \cB/c_0(\cGh)$ is denoted by $\om_\infty$.
\end{proposition}
\begin{proof}
It suffices to observe that $\psi_{x+y} \psi_{x+y,x} = \psi_x$ for all $x,y \in \Irred(\cG)$.
\end{proof}

Denote by $(\cB_\infty,\om_\infty)\dpr$ the von Neumann algebra
generated by $\cB_\infty$ in the GNS-construction for the state
$\om_\infty$. It is easy to check that $\om_\infty$ is a KMS state on
$\cB_\infty$ with modular group given by
$\sigma_t^{\om_\infty}(\psi_{\infty,x}(A)) = \psi_{\infty,x}(Q_x^{it}
A Q_x^{-it})$ for all $t \in \R$, $x \in \Irred(\cG)$ and $A \in
\B(H_x)$. In particular, $\om_\infty$ induces a normal faithful state
on $(\cB_\infty,\om_\infty)\dpr$.

\begin{theorem}
Denote $\cG=A_o(F)$ and suppose $\cG \not\cong \SU(2),\SU_{-1}(2)$. Let $\mu$ be any generating measure on $\Irred(\cG)$. The linear map
$$T : \cB_\infty \recht \ell^\infty(\cGh) : T(a) = (\id \ot \om_\infty)\be_\infty(a)$$
yields a $^*$-homomorphism $T : \cB_\infty \recht H^\infty(\cGh,\mu)$ satisfying $\epsh T = \om_\infty$, and a $^*$-isomorphism
$$(\cB_\infty,\om_\infty)\dpr \cong (H^\infty(\cGh,\mu),\epsh) \; .$$
\end{theorem}

\begin{proof}
We first claim that $(\psi_x \ot \om_\infty)\be_\infty(a) = \om_\infty(a)$ for all $a \in \cB_\infty$. It then follows that $T(a) \in
H^\infty(\cGh,\mu)$ for all $a \in \cB_\infty$. To prove our claim, observe that, for all $x$ and $y \geq x$, $(\psi_x \ot \psi_y)\deh$ is a convex
combination of $\psi_{y-x},\ldots,\psi_{y+x}$. It follows that $(\psi_x \ot \om)\be(a) = \om(a)$ for all $a \in \cB$.

To show that $T$ is multiplicative, it suffices to show that, for all fixed $a,b \in \cB_\infty$,
\begin{equation}\label{eq.multiplicative}
\| ( T(a)T(b) - T(ab))p_n \| \recht 0
\end{equation}
when $n \recht \infty$. Indeed, for fixed $a,b \in \cB_\infty$ and a fixed $x$,
$$(T(a) \cdot T(b))p_x = \slim_{n \recht \infty} (\id \ot \psi_\mu^{\star n})\deh(T(a)T(b))p_x \; .$$
By the transience of the state $\psi_\mu$, the expression on the right, for $n$ big, only takes into account $T(a)T(b)p_m$ for $m$ big. This last
expression is close to $T(ab)p_m$. But, $(\id \ot \psi_\mu^{\star n})\deh(T(ab)) = T(ab)$ and we are done.

Hence, to prove the multiplicativity of $T$, it remains to show \eqref{eq.multiplicative}. It suffices to show that for all $x$ and $A \in \B(H_x)$
fixed,
\begin{equation}\label{eq.sufficient}
\| (\id \ot \om_\infty)\be_\infty(\psi_{\infty,x}(A))p_n - \psi_{n,x}(A) \| \recht 0
\end{equation}
when $n \recht \infty$. Fix $x$ and $A \in \B(H_x)$. Take $y \geq x$ and $z \geq y$. Then,
\begin{align*}
(\id \ot \psi_{x+z}) & \deh(\psi_{\infty,x}(A)) p_{y} = \sum_{s \in y \ot z} (\id \ot \psi_{x+z})\bigl( V(y \ot (x+z),x+s) \psi_{x+s,x}(A) V(y \ot (x+z),x+s)^* \bigr) \\
&= \sum_{s \in y \ot z} \frac{[x+s+1]}{[y+1] [x+z+1]} V((x+s) \ot (x+z),y)^* (\psi_{x+s,x}(A) \ot 1) V((x+s) \ot (x+z),y) \; .
\end{align*}
From Lemma \ref{lem.crucial-other-side}, we get a constant $C$, only depending on $q$, such that
$$d_\T \bigl( (V(x \ot s,x+s) \ot 1) V((x+s) \ot (x+z),y), (1 \ot V(s \ot (x+z),y-x))V(x \ot (y-x),y) \bigr) \leq C q^{-x-z+s} \; .$$
Note however that this statement only makes sense when $y-x \in s \ot (x+z)$. This is the case for $s \geq z-y+2x$ and so, we can safely go on
because our estimate is bigger than $1$ if $s < z-y+2x$. Hence, we get a constant $D$ such that
$$(\id \ot \psi_{x+z})\deh(\psi_{\infty,x}(A)) p_{y} \approx \sum_{s \in y \ot z} \frac{[x+s+1]}{[y+1] [x+z+1]} \psi_{y,x}(A) = \psi_{y,x}(A)$$
with error $\leq \sum_{s \in y \ot z} 2\|A\| C \frac{[x+s+1] q^{x+s}}{[y+1] [x+z+1] q^{x+z}} q^{-x} \leq Dq^{-x} \|A\| \frac{2y+1}{[y+1]}$. Since
this estimate holds for all $z \geq y$, we find that
$$\bigl\| (\id \ot \om_\infty)\be_\infty(\psi_{\infty,x}(A)) p_y - \psi_{y,x}(A) \| \leq D q^{-x} \|A\| \frac{2y+1}{[y+1]} \; .$$
So, \eqref{eq.sufficient} follows and the multiplicativity of $T$ has been proven.

It is obvious that $\epsh T = \om_\infty$.
Consider the adjoint action of $\cG$ on $\ell^\infty(\cGh)$ given by
$\Phi(a) = \Vb(a \ot 1)\Vb^*$ for all $a \in \ell^\infty(\cGh)$. Since
$$U^{x+y} (\psi_{x+y,x}(A) \ot 1) (U^{x+y})^* = (\psi_{x+y,x} \ot \id)(U^x (A \ot 1)(U^x)^*) \; ,$$
the action $\Phi$ restricts to an action of $\cG$ on $\cB$. Moreover, the action $\Phi$ preserves the ideal $c_0(\cGh) \subset \cB$, yielding an
action $\Phi_\infty$ of $\cG$ on $\cB_\infty$. We have $(\id \ot h)\Phi_\infty(a) = \om_\infty(a)1$ for all $a \in \cB_\infty$. So, $\Phi_\infty$ is
an ergodic action and $\om_\infty$ is the unique invariant state. By definition,
$$\Phi_\infty(\psi_{\infty,x}(A)) = (\psi_{\infty,x} \ot \id)(U^x (A
\ot 1)(U^x)^*) \; .$$ From Lemma \ref{lem.injective-psi}, we know that $\psi_{\infty,x} : \B(H_x) \recht \cB_\infty$ is an injective linear map. So,
we conclude that the irreducible representation $U^x$ appears with multiplicity one in $\Phi_\infty$ when $x$ is even and with multiplicity zero when
$x$ is odd. Moreover, the $^*$-homomorphism $T$ intertwines $\Phi_\infty$ with the adjoint action of $\cG$ on $H^\infty(\cGh,\mu)$.

From Corollary 3.5 in \cite{INT}, we know that the multiplicities of the irreducible representations in the adjoint action on $H^\infty(\cGh,\mu)$
are at most the multiplicities in $\Phi_\infty$. Since $\om_\infty$ yields a faithful, normal state on $(\cB_\infty,\om_\infty)\dpr$ and since $\epsh T = \om_\infty$, the homomorphism $T : (\cB_\infty,\om_\infty)\dpr \recht H^\infty(\cGh,\mu)$ is faithful. It follows
that $T$ is a $^*$-isomorphism.
\end{proof}

\subsection*{Martin boundary}

The Martin boundary and the Martin compactification of a discrete
quantum group have been defined by Neshveyev and Tuset in
\cite{NT}. We first introduce the necessary terminology and notation
and then prove that the Martin compactification of the dual of
$A_o(F)$ is equal to the compactification $\cB$ constructed above.

Let $\cGh$ be a discrete quantum group and
let $\mu$ be a probability measure on $\Irred(\cG)$. We have an
associated Markov operator $Q_\mu$ and a classical random walk on $\Irred(\cG)$ with $n$-step transition probabilities
given by \eqref{eq.nstep}. We suppose throughout that $\mu$ is a
generating measure and that $\mu$ is transient.
It follows that $0 < \sum_{n=1}^\infty p_n(x,y) < \infty$ for all $x,y \in
\Irred(\cG)$.

Denote by $c_c(\cGh) \subset c_0(\cGh)$ the algebraic direct sum of
the algebras $\B(H_x)$. We define, for $a \in c_c(\cGh)$,
$$G_\mu(a) = \sum_{n=0}^\infty Q_\mu^n(a) \; .$$
Observe that usually $G_\mu(a)$ is unbounded, but it makes sense in
the multiplier algebra of $c_c(\cGh)$, i.e.\ $G_\mu(a)p_x \in \B(H_x)$
makes sense for every $x \in \Irred(\cG)$. Moreover, $G_\mu(p_\eps)$ is strictly positive and
central. This allows to define the Martin kernel as follows.

Whenever $\mu$ is a measure on $\Irred(\cG)$, we use the notation
$\omu$ to denote the measure given by $\omu(x) = \mu(\ox)$.

\begin{definition}[Defs.\ 3.1 and 3.2 in \cite{NT}]
Define
$$K_\mu : c_c(\cGh) \recht \ell^\infty(\cGh) : K_\mu(a) = G_\mu(a)
G_\mu(p_\eps)^{-1} \; .$$
Define the \emph{Martin compactification} $\Atil_\mu$ as the
\cst-subalgebra of $\ell^\infty(\cGh)$ generated by
$K_{\omu}(c_c(\cGh))$ and $c_0(\cGh)$. Define the \emph{Martin
  boundary} $A_\mu$ as the quotient $\Atil_\mu / c_0(\cGh)$.
\end{definition}

The aim of this section is to prove the following result.

\begin{theorem}
Denote $\cG=A_o(F)$ and suppose $\cG \not\cong
\SU(2),\SU_{-1}(2)$. Let $\mu$ be a generating measure on
$\Irred(\cG)$ with finite first moment~:
$$\sum_{x \in \N} x \mu(x) < \infty \; .$$
Then, the Martin compactification $\Atil_\mu$ equals the
compactification $\cB$ defined in Proposition \ref{prop.compactif}. In
particular, the Martin boundary $A_\mu$ equals $\cB_\infty$.
\end{theorem}

\begin{proof}
Introduce the notation $p_x g_\mu(x,y) = p_x G_\mu(p_y)$. One has
$g_\mu(0,x) = g_{\omu}(x,0) \dimq(x)^2$. So, if $A \in \B(H_x)$ and $y
\geq x$, we get
$$G_{\omu}(Ap_x) p_y = \sum_{z \in x \ot y} g_{\omu}(0,z) (\id \ot
\psi_z)(V(y \ot z,x) A V(y \ot z,x)^*) \; .$$ An easy computation yields
\begin{equation*}
K_{\omu}(Ap_x) p_y = \sum_{z = 0}^x \frac{g_\mu(y-x+2z,0)}{g_\mu(y,0)} \frac{\dimq(y-x+2z) \dimq(x)}{\dimq(y)} V(x \ot (y-x+2z))^*(A \ot 1) V(x \ot
(y-x+2z)) \; .
\end{equation*}
From Proposition 4.7 in \cite{NT}, we know that
$$\lim_{x \recht \infty} \frac{g_\mu(x+1,0)}{g_\mu(x,0)} = q^2 \; .$$
Using the Notation \ref{not.psipsi} below, the previous formula and
the asymptotics for the quantum dimensions, we find that, for all $x
\in \Irred(\cG)$ and $A \in \B(H_x)$,
$$\lim_{y \recht \infty} \|K_{\omu}(Ap_x)p_y - \sum_{z=0}^x q^{-x+2z}
\dimq(x) \psi_{y,x}^z(A) \| = 0 \; .$$ Denoting by $[Y]$ the closed linear span of $Y$, we conclude that
\begin{equation} \label{eq.wezijnerbijna}
[c_0(\cGh) + K_{\omu}(c_c(\cGh))] = [c_0(\cGh) + \sum_{z=0}^x
q^{2z} \psi_{\infty,x}^z(A) \mid x \in \N, A \in \B(H_x)] \; .
\end{equation}
Recall the \cst-algebra $\cB$ defined in Proposition
\ref{prop.compactif}. Combining \eqref{eq.wezijnerbijna} with
\eqref{eq.hophop} in Lemma \ref{lem.heeelp} below, we conclude that
$[c_0(\cGh) + K_{\omu}(c_c(\cGh))] \subset \cB$. The opposite
inclusion follows by combining \eqref{eq.wezijnerbijna} with
\eqref{eq.hiphip} in Lemma \ref{lem.heeelp} below. In particular, the
\cst-algebra generated by $c_0(\cGh)$ and $K_{\omu}(c_c(\cGh))$ equals $\cB$.
\end{proof}

\begin{notation} \label{not.psipsi}
Definition \ref{def.psi} admits the following natural
generalization. For all $y \geq x \geq z$, we write
$$\psi_{y,x}^z : \B(H_x) \recht \B(H_y) : \psi_{y,x}^z(A) = V(x \ot
(y-x + 2z),y)^* (A \ot 1) V(x \ot (y-x + 2z),y) \; .$$ Note that $\psi_{y,x} = \psi_{y,x}^0$. We write as well $\psi_{\infty,x}^z(A)$ defined for $x
\geq z$ and $A \in \B(H_x)$ by
$$\psi_{\infty,x}^z(A) p_y = \psi_{y,x}^z(A) \quad\text{whenever}\quad
y \geq x \; .$$
\end{notation}

\begin{lemma} \label{lem.heeelp}
There exists a constant $C > 0$, only depending on $q$ such that, for all
$x,y,z \in \N$ and $A \in \B(H_x)$, we have
\begin{align}
\| \psi_{x+y+z,x+y} \psi_{x+y,x}^r (A) - \psi_{x+y+z,x}^r(A) \| & \leq C
q^{y+r} \|A\| \quad\text{for all}\; 0 \leq r \leq x \; ,
\label{eq.hophop} \\
\| \psi_{x+y+z,x+y}^r \psi_{x+y,x}(A) - \psi_{x+y+z,x}(A) \| & \leq C
q^{y-r} \|A\| \quad\text{for all}\; 0 \leq r \leq x+y \; . \label{eq.hiphip}
\end{align}
\end{lemma}

\begin{proof}
Inequality \eqref{eq.hophop} follows from \eqref{eq.best2} in Lemma
\ref{lem.crucial-improved}, while \eqref{eq.hiphip} follows from
\eqref{eq.best-a2} in Lemma \ref{lem.crucial-other-side}.
\end{proof}

\section{A general exactness result} \label{sec.generalexact}

In \cite{BDRV} a notion of monoidal equivalence of compact quantum groups was introduced. It was shown in particular that for all $F \in \GL(n,F)$
with $F \overline{F} = c 1$, $c = \pm 1$, $\Tr(F^*F) = q + \frac{1}{q}$ and $0 < q < 1$, the quantum groups $A_o(F)$ and $\SU_{-cq}(2)$ are
monoidally equivalent.

In this section we prove that the exactness of the reduced \cst-algebra of a compact quantum group is invariant under monoidal equivalence. As a
corollary we obtain an alternative proof for the first half of Corollary \ref{cor.exactAO}.

\begin{theorem}
Let $\cG$ and $\cG_1$ be compact quantum groups. Let $\vphi : \cG \recht \cG_1$ be a monoidal equivalence in the sense of \cite{BDRV}, Definition
3.1. Let $B\red$ be the associated reduced \cst-algebra. The following statements are equivalent.
\begin{itemize}
\item $C(\cG)\red$ is exact.
\item $C(\cG_1)\red$ is exact.
\item $B\red$ is exact.
\end{itemize}
\end{theorem}

\begin{proof}
Following \cite{BDRV}, Theorem 3.9, we consider the $^*$-algebra $\cB$ generated by the coefficients of unitary elements $X^x \in
\B(H_x,H_{\vphi(x)}) \ot \cB$. We consider the canonical invariant state $\om$ on $\cB$ and denote by $B\red$ the associated reduced \cst-algebra.

By symmetry, i.e.\ using the inverse $\vphi^{-1} : \cG_1 \recht \cG$, we consider the $^*$-algebra $\cBtil$ generated by the coefficients of unitary
elements $Y^x \in \B(H_{\vphi(x)},H_x) \ot \cBtil$. We denote by $\omtil$ the invariant state on $\cBtil$ and by $\Btil\red$ the associated reduced
\cst-algebra. Observe that we have a canonical anti-isomorphism $\pi : B\red \recht \Btil\red$ given by $(\id \ot \pi)(X^x) = (Y^x)^*$ for all $x \in
\Irred(\cG)$. So, $\Btil\red$ is nothing else than the opposite of $B\red$.

Suppose first that $B\red$ is exact. Then, also $\Btil\red$ is exact. Moreover, we get an injective $^*$-homomorphism $\theta : C(\cG)\red \recht
\Btil\red \ot B\red$ given by $(\id \ot \theta)(U^x) = Y^x_{12} X^x_{13}$. A priori, $\theta$ defines a $^*$-homomorphism $C(\cG)_{\text{\rm u}}
\recht \Btil\red \ot B\red$, but it is easy to verify that $(\omtil \ot \om)\theta$ is the Haar state. So, $\theta$ is well defined on $C(\cG)\red$.
It follows that $C(\cG)\red$ is exact. In a similar way, we deduce that $C(\cG_1)\red$ is exact.

Suppose next that $C(\cG)\red$ is exact. Then, $\cGh$ is an exact quantum group. Since $B\red$ is Morita equivalent with a reduced crossed product
$\cGh \ltimesred \cK$, it follows that $B\red$ is exact. In a similar way, exactness of $C(\cG_1)\red$ implies exactness of $B\red$.
\end{proof}

\begin{corollary}
Let $\cG = A_o(F)$. Then, $C(\cG)\red$ is exact.
\end{corollary}
\begin{proof}
By amenability of the (dual of) $\SU_q(2)$, the exactness of its reduced (=universal) \cst-algebra is obvious. The result follows since every
$A_o(F)$ is monoidally equivalent with some $\SU_q(2)$.
\end{proof}

The same argument admits the following generalization.

\begin{corollary}
The reduced \cst-algebra of any compact quantum group monoidally equivalent with a $q$-deformation of a simple compact Lie group, is exact.
\end{corollary}

\begin{remark}
We can as well give a sledgehammer argument for the exactness of the reduced \cst-algebra of $A_u(F)$. It follows from \cite{BDRV} that any $A_u(F)$
is monoidally equivalent with an $A_u(F)$ with $F \in \GL(2,\C)$. But, it follows from \cite{banica2} that the reduced \cst-algebra of such an
$A_u(F)$ is a subalgebra of the reduced free product of $\SU_q(2)$ and $S^1$ and hence, we are done.
\end{remark}

\section{Factoriality and simplicity} \label{sec.factoriality}

We prove that, at least in most cases, the von Neumann algebras associated with $\cG = A_o(F)$ are factors. We determine their Connes invariants and
we prove that the reduced \cst-algebras are simple. In combination with the results above on the Akemann-Ostrand property, we obtain new examples of
generalized solid, in particular prime, factors.

On the von Neumann algebra side, we get
\begin{theorem} \label{thm.mainvnalg}
Let $N \geq 3$, $F \in \GL(N,\C)$ with $F \overline{F} = \pm 1$. Suppose that $\|F\|^2 \leq \Tr(FF^*)/\sqrt{5}$. Write $\cG = A_o(F)$ and $M =
C(\cG)\vnalg$.
\begin{itemize}
\item $M$ is a full generalized solid factor with almost periodic state $h$. In particular, $M$ is prime.
\item $\Sd(M)$ is the subgroup $\Gamma$ of $\R^*_+$ generated by the
  eigenvalues of $Q \ot Q^{-1}$. In particular, $M$ is of type II$_1$ when $FF^* = 1$~; of
type III$_\lambda$ when $\Gamma = \lambda^\Z$ and of type III$_1$ in the other cases.
\item The \cst-subalgebra of $\B(L^2(\cG))$ generated by $\lambda(C(\cG))$ and $\rho(C(\cG))$ contains the compact operators.
\end{itemize}
If $FF^* = 1$ and $N \geq 3$, $M$ is a solid, in particular prime, II$_1$ factor.
\end{theorem}

On the C$^*$-algebra side, we obtain
\begin{theorem} \label{thm.simplicity}
Let $N \geq 3$, $F \in \GL(N,\C)$ with $F \overline{F} = \pm 1$. Suppose that $\|F\|^8 \leq \frac{3}{8} \Tr(FF^*)$. Write $\cG = A_o(F)$. Then
$C(\cG)\red$ is a simple exact \cst-algebra and $h$ is the unique state on $C(\cG)\red$ satisfying the KMS condition with respect to $(\sigma_t^h)$.
In particular, for $\cG=A_o(I_N)$ and $N \geq 3$, $C(\cG)\red$ is a simple exact \cst-algebra with unique tracial state $h$.
\end{theorem}

Theorems \ref{thm.mainvnalg} and \ref{thm.simplicity} are proven through a careful analysis of the quantum analogue of the operation of \lq
conjugation by the generators\rq\ in free groups, see Definition \ref{def.operators}.

Fix a matrix $F \in \GL(n,\C)$ satisfying $F \overline{F} = \pm 1$ and put $\cG = A_o(F)$. Define $H_1 = \C^n$ and $U^1:=U$, the fundamental
representation on $H_1$. Recall that the modular theory of the compact quantum group $\cG$ is encoded by positive invertible elements $Q_x \in
\B(H_x)$. In the case of $\cG= A_o(F)$, we write $Q:=Q_1$ and we have
$$Q = F^t \overline{F} \quad\text{and}\quad Q^{-1} = F F^* \; .$$
Write $F \overline{F} = c1$ with $c = \pm 1$. Write $t_1 = \Tr(Q)^{-1/2} \sum_{i=1}^n e_i \ot Fe_i$, which is a \emph{unit} invariant vector for the
tensor square $U^{\ot 2}$.

In order to study factoriality and simplicity, we introduce the following operators, using Notation \ref{not.GNSdiscrete}. Recall as well the regular
representation $\rho : C(\cG) \recht \B(L^2(\cG))$ given by \eqref{eq.regrho}. We denote the anti-homomorphism $\rhoop$ defined by
$$\rhoop(a)\rho(b)\xi_0 = \rho(ba) \xi_0 \quad\text{for all}\quad a \in \Calg(\cG), b \in C(\cG) \; ,$$
where $\Calg(\cG) \subset C(\cG)$ is the dense $^*$-subalgebra given by the coefficients of finite-dimensional representations of $\cG$. Note that
$\rhoop$ is not involutive~: we have $\rhoop(a)^* = \rhoop(\sigma_i^h(a)^*)$ where the modular group $(\sigma_t^h)$ is given by
$$(\id\otimes\sigma_t^h)(U) = (Q^{it}\otimes 1)U(Q^{it}\otimes 1) \; .$$

\begin{definition} \label{def.operators}
We define operators $T$ and $\Ptil$ as follows.
\begin{align*}
T &: L^2(\cG) \recht H_1 \ot L^2(\cG) \ot H_1 : \text{flip} \circ T = \frac{1}{\dimq(1)} \bigl((\Lambdahl\otimes\rho)(U) -
(\Lambdahl\otimes\rhoop)(U)\bigr) \; , \\
\Ptil &: C(\cG) \recht C(\cG) : \Ptil(a) = \frac{1}{2 \Tr(Q)} \bigl((\Tr_Q \ot \id)(U^*(1 \ot a)U) + (\Tr_{Q^{-1}} \ot \id)(U(1 \ot a)U^*)\bigr) \; ,
\end{align*}
where $\text{flip}$ denotes the identification
$$\text{flip} : H_1 \ot H_1 \ot L^2(\cG) \recht H_1 \ot L^2(\cG) \ot
H_1 : \text{flip}(\xi \ot \eta \ot \mu) = \eta \ot \mu \ot \xi \; .$$ We also use $\Ptil$ on the Hilbert space level, writing
$$P : L^2(\cG) \recht L^2(\cG) : P \rho(a)\xi_0 = \rho(\Ptil(a)) \xi_0 \quad\text{for all}\quad a \in C(\cG) \; .$$
\end{definition}

It is straightforward to check that
\begin{align*}
T &= \frac{c}{\dimq(1)^{1/2}} \sum_{i,j=1}^n Fe_j \ot (\rho(U_{ij}) - \rhoop(U_{ij})) \ot e_i \; , \\
T^*T &= 2(1- P) \; .
\end{align*}

\begin{remark}
The relevance of the operator $T$ in the study of the factoriality of $C(\cG)\vnalg$ is clear. Indeed, $C(\cG)\vnalg$ is a factor if and only if $T
\eta = 0$ implies $\eta \in \C \xi_0$.
\end{remark}

Together with proving the factoriality of $C(\cG)\vnalg$, we compute Connes' invariants for $C(\cG)\vnalg$. In order to do so, we introduce the
following deformation of $T$.
$$T_t = \frac{1}{\dimq(1)}\bigl((\Lambdahl\ot\rho\sigma_t^h)(U) - (\Lambdahl\ot\rhoop)(U)\bigr) =
\frac{c}{\dimq(1)^{1/2}} \sum_{i,j=1}^n Fe_j \ot (\rho\sigma_t^h(U_{ij}) - \rhoop(U_{ij})) \ot e_i \;$$

The following is the major technical result of the section. It follows from a series of lemmas proved at the end of the section. We already deduce
factoriality and simplicity results from it.

\begin{proposition} \label{prop.mainestimate}
If $\|Q\| \leq \frac{\Tr(Q)}{\sqrt{5}}$, there exist $C_1 > C_2 > 0$ such that
$$\|T_s \xi \| \geq \sqrt{C_1^2 \|\xibol\|^2 + D_s
|\la \xi_0,\xi \ra|^2} - C_2 \|\xibol\| \; ,$$ for all $s \in \R$, $\xi \in L^2(\cG)$, where
$$\xi = \xibol + \xi_0 \la \xi_0,\xi \ra \quad\text{and}\quad D_s =
2(1-|\la(1 \ot Q^{is})t,t\ra|^2) \; .$$
\end{proposition}

We already show how theorems \ref{thm.mainvnalg} and \ref{thm.simplicity} follow from Proposition \ref{prop.mainestimate}.

\begin{proof}[Proof of Theorem \ref{thm.mainvnalg}]
Write $S = C(\cG)\red$ and consider the operator $P \in \B(L^2(\cG)$ introduced in Definition \ref{def.operators}. Note that $2(1-P) = T^*T$. From
the definition of $P$, we get that $P \in C^*(\lambda(S),\rho(S))$. From Proposition \ref{prop.mainestimate}, we get $0 < C < 1$ such that the
spectrum of $P$ is included in $[0,C] \cup \{1\}$ and the spectral projection of $\{1\}$ is precisely the projection onto $\C \xi_0$. It follows that
$C^*(\lambda(S),\rho(S))$ contains the compact operators and that $M = C(\cG)\vnalg$ is a full factor. From Corollary \ref{cor.solid}, we already
know that $M$ is a generalized solid von Neumann algebra. Combining with Proposition \ref{prop.solidprime}, we get that $M$ is prime.

Denote by $\Gamma$ the subgroup of $\R^*_+$ generated by the eigenvalues of $Q \ot Q^{-1}$. In order to show that $\Sd(M) = \Gamma$ it suffices to
show that, given a sequence $(s_n)$ in $\R$, $\sigma_{s_n}^h \recht 1$ in $\Out(M) = \frac{\Aut(M)}{\Inn(M)}$ if and only if $$|\la (1 \ot
Q^{is_n})t,t \ra| \recht 1 \; .$$ One implication being obvious, suppose that $\sigma_{s_n}^h \recht 1$ in $\Out(M)$. Take unitaries $u_n \in M$ such
that $(\Ad u_n) \sigma_{s_n}^h \recht \id$ in $\Aut(M)$. It follows that
$$\| T_{s_n} \rho(u_n^*) \xi_0 \| \recht 0 \; .$$
Applying Proposition \ref{prop.mainestimate}, we first get that $\|(\rho(u_n^*)\xi_0)^\circ\| \recht 0$ and next that $|\la (1 \ot Q^{is_n})t,t \ra|
\recht 1$.
\end{proof}

\begin{proof}[Proof of Theorem \ref{thm.simplicity}]
Consider the operator $\Ptil : C(\cG)\red \recht C(\cG)\red$ as in Definition \ref{def.operators}. From Proposition \ref{prop.mainestimate}, we get a
constant $0 < C < 1$ such that $\|\Ptil(a)\|_2 \leq C \|a\|_2$ for all $a \in C(\cG)\red$ with $h(a) = 0$. If $C(\cG)_n$ denotes the linear span of
matrix coefficients of $U^0$, ..., $U^n$, we have $\Ptil^k(C(\cG)_n) \subset C(\cG)_{n+2k}$. In particular for $a\in C(\cG)_n$ such that $h(a)=0$ we
have
$$
\|\Ptil^k(a)\| \leq p(n+2k) \|Q\|^{n+2k} \|\Ptil^k(a)\|_2 \leq p(n+2k) \|Q\|^n (C\|Q\|^2)^{k} \|a\|_2 \; ,
$$
where $p$ is a fixed polynomial given by the Property of Rapid Decay for $A_o(F)$ (see Remark \ref{rem.rapiddecay} below). Hence if $C\|Q\|^2 < 1$ we
find that $\Ptil^k(a) \recht 0$. Since $\Ptil$ is unital, this implies that $\Ptil^k(a) \to h(a)$ for any $a \in C(\cG)\red$, hence $a$ cannot be in
a non-trivial ideal. Moreover $\Ptil$ leaves invariant any state $\varphi$ satisfying the KMS condition with respect to $(\sigma_t^h)$, hence $h$ and
$\varphi$ agree on any $a \in C(\cG)\red$. It follows from \eqref{eq.daargaathetom} in Remark \ref{rem.goedgeteld} that the condition $C\|Q\|^2 < 1$
is satisfied whenever $\|Q\|^4 \leq \frac{3}{8} \Tr(Q)$.
\end{proof}

\begin{remark} \label{rem.rapiddecay}
In the proof of the previous theorem, we made use of the Property of Rapid Decay for universal quantum groups as introduced in \cite{ver2}. This
property yields a control over the norm in $C(\cG)\red$ using the norm in $L^2(\cG)$.

Denote by $C(\cG)_n$ the linear span of matrix coefficients of $U^0$, ..., $U^n$, and by $\|\,\cdot\,\|_2$ the GNS norm associated with $h$. One
possible definition of Property RD goes as follows:
$$
  \exists p\in\R[X] \quad\text{such that}\quad \forall n\in\N \; ,
  \;\; a\in C(\cG)_n \; , \;\; \|\rho(a)\| \leq p(n) \|a\|_2 \; .
$$
In the case where $\cG = A_o(F)$, it is proven in Theorem 3.9 in \cite{ver2} that property RD holds if and only if $F$ is a multiple of a unitary
matrix, i.e.\ $Q=1$. In fact the techniques of \cite{ver2} still work in the non-unimodular case, but yield non-polynomial bounds: we get a
polynomial $p \in \R[X]$ such that
$$
\|\rho(a)\| \leq
  \|Q\|^n p(n) \|a\|_2 \; ,
$$
for all $n \in \N$ and all $a \in C(\cG)_n$.
\end{remark}

So, it remains to prove Proposition \ref{prop.mainestimate}, which takes the rest of the section.

The proof is not very hard, but somehow computationally involved. In order to streamline our computations, we choose explicit representatives for the
irreducible representations of $\cG = A_o(F)$, as well as for the intertwiners $V(x \ot y, z)$ with tensor products of irreducible representations.
We know that $\Mor(1^{\ot n},1^{\ot n})$ is isomorphic with the Temperley-Lieb algebra. In particular, we have the Jones-Wenzl projection $p_n \in
\Mor(1^{\ot n},1^{\ot n})$, which allows to define $H_n := p_n H_1^{\ot n}$ and take $U^n$ as the restriction of the $n$-fold tensor product $U^{\ot
n}$ to $H_n$. We write $1_n = 1^{\ot n}$.

Using Theorem 3.7.1 in \cite{SFS}, we can recursively define the \emph{unit} vectors $t_x \in \Mor(x \ot x,0)$ and the \emph{isometries} $V((x+z) \ot
(z+y),x+y) \in \Mor((x +z) \ot (z+y),x+y)$ using the formulas
\begin{align}
t_1 &= \Tr(F^*F)^{-1/2} \sum_{i = 1}^n e_i \ot Fe_i \; , \\
t_{x+y} &= \Bigl(\frac{[x+1] [y+1]}{[x+y+1]} \Bigr)^{1/2} (p_{x+y} \ot
p_{y+x})(1_x \ot t_y \ot 1_x)t_x  \; , \label{eq.relative} \\
V((x+z) \ot (z+y),x+y) &= \Bigl(\frac{[z+1] \qbinom{x+z}{z}
  \qbinom{y+z}{z}}{\qbinom{x+y+z+1}{r}}\Bigr)^{1/2} (p_{x+z} \ot
p_{z+y})(1_x \ot t_z \ot 1_y)p_{x+y} \; . \label{eq.innerproduct}
\end{align}

Here we used the usual notation of $q$-numbers, $q$-factorials and $q$-binomial coefficients. We take as before $0 < q < 1$ such that
$$\Tr(F^*F) = q + \frac{1}{q} \; , \quad
F \overline{F} = c 1 \quad\text{where}\quad c = \pm 1 \; .
$$
Then, we use the following notation.

\begin{notation}
With $0 < q < 1$ fixed, we write the $q$-numbers, $q$-factorials and $q$-binomial coefficients.
$$
[n] = \frac{q^{n} - q^{-n}}{q - q^{-1}} \; , \quad [n]! = [n] [n-1] \cdots [1] \; , \quad \qbinom{n}{r} = \frac{[n]!}{[r]! [n-r]!} \; .$$
\end{notation}

Note that $\dimq(x) = [x+1]$. Wenzl's recursion formula for the projections $p_n$ admits the following generalisation (see equation (3.8) in
\cite{FK}, page 462).
\begin{equation}\label{eq.recursion}
p_n = \Bigl( 1 - \sum_{k=1}^{n-1} (-c)^{n-1-k} \frac{[2] [k]}{[n]} (1_{k-1} \ot t_1 \ot 1_{n-1-k} \ot t_1^*) \Bigr)(p_{n-1} \ot 1) \; .
\end{equation}
Note that, multiplying on the left with $p_{n-1} \ot 1$, we obtain Wenzl's recursion
\begin{equation}\label{eq.wenzl}
p_n = p_{n-1} \ot 1 - \frac{[2] [n-1]}{[n]}(p_{n-1} \ot 1)(1_{
  n-2} \ot t_1t_1^*)(p_{n-1} \ot 1) \; .
\end{equation}

\begin{notation} \label{not.factintertwin}
The study of $T$ consists in comparing the left and right actions of the coefficients of $U$, and hence it has a natural counterpart at the level of
representations. More precisely, let us introduce the following short hand notations.
\begin{alignat*}{2}
\phipl &:= V(1 \ot (x+1),x) \; , \qquad
& \phipr &:= V((x+1) \ot 1,x) \; , \\
\phiml &:= V(1 \ot (x-1),x) \; , \qquad & \phimr &:= V((x-1) \ot 1,x) \; .
\end{alignat*}
For any $x\in\N$, we also define $\sigma : H_x\ot H_1 \recht H_1\ot H_x$ by $\sigma (\eta\ot\mu) = \mu\ot\eta$.
\end{notation}

\begin{lemma} \label{lem.identifications}
We have $T = T^+ + T^-$, where $T^+(H_x \ot H_x) \subset H_{x+1} \ot H_{x+1}$ and $T^-(H_x \ot H_x) \subset H_{x-1} \ot H_{x-1}$, are defined by the
formulas
\begin{align*}
T^+ \eta &= \sqrt{\frac{\dimq(x+1)}{\dimq(1)\dimq(x)}} \;\; \bigl((Q^{-1} \ot 1) \phipl \ot \phipr \; - \;  \si\phipr \ot \si^*(Q^{-1} \ot
1)\phipl \bigr) \eta  \; , \\
T^- \eta &= \sqrt{\frac{\dimq(x-1)}{\dimq(1)\dimq(x)}} \;\; \bigl((Q^{-1} \ot 1) \phiml \ot \phimr  \; - \; \si\phimr \ot \si^*(Q^{-1} \ot 1)\phiml
\bigr) \eta \; ,
\end{align*}
for all $\eta \in H_x \ot H_x$. The corresponding formulas for $T_s$ are obtained by composing the first term of each difference with $(Q^{is} \ot 1)
\ot (1 \ot Q^{-is})$.
\end{lemma}

\begin{proof}
By definition of the tensor product of representations of $\cG$, we have
\begin{align*}
(\omega_{\eta,\xi}\ot\rho)(U^x)(\omega_{\eta',\xi'}\ot\rho)(U^{x'}) &= \sum_{y \in x\ot x'} (\omega_{V(x\ot x',y)^*(\eta\ot\eta'), V(x\ot
x',y)^*(\xi\ot\xi')}\ot\rho)(U^y) \; .
\end{align*}
Since the definition of $T$ involves multiplication on the left and multiplication on the right by coefficients of $U^1$ and since $1 \ot x$ and $x
\ot 1$ split into a direct sum of $x-1$ and $x+1$, clearly $T$ consists of four terms, say $T = T^+ + T^-$ and $T^\pm = T^\pm_l - T^\pm_r$. More
explicitly,
\begin{equation}\label{eq.daargaanwe}
T^+_l \rho\bigl((\om_{\eta,\xi} \ot \id)(U^x)\bigr) \xi_0 = \frac{c}{\dimq(1)^{1/2}} \sum_{i,j=1}^n Fe_j \ot \rho\bigl( (\om_{(\phiml)^*(e_i \ot
\eta),(\phiml)^*(e_j \ot \xi)} \ot \id)(U^{x+1}) \bigr) \xi_0 \ot e_i \; ,
\end{equation}
and $T^+_r$, $T^-_l$ and $T^-_r$ are defined analogously. From \eqref{eq.regrho}, we get that
$$\rho\bigl((\om_{\eta,\xi} \ot \id)(U^x)\bigr) \xi_0 = \xi \ot (1 \ot
\eta^*)t_x \; .$$ To compute the right hand side of \eqref{eq.daargaanwe}, observe that
\begin{equation} \label{eq.eerstedeel}
\sum_{i=1}^n (1 \ot (e_i \ot \eta)^* \phiml)t_{x+1} \ot e_i = (1 \ot 1 \ot \eta^*)(1 \ot \phiml)t_{x+1} = (1 \ot 1 \ot \eta^*) (\phipr \ot 1)t_x =
\phipr (1 \ot \eta^*)t_x \; .
\end{equation}
Using the equality $\sum_{j=1}^n Fe_j \ot e_j = c \dimq(1)^{1/2} (Q^{-1} \ot 1)t_1$, we observe as well that
\begin{equation} \label{eq.tweededeel}
\sum_{j=1}^n Fe_j \ot (\phiml)^* (e_j \ot \xi) = c \dimq(1)^{1/2} (Q^{-1} \ot (\phiml)^*)(t_1 \ot \xi) = c \sqrt{\frac{\dimq(x+1)}{\dimq(x)}} (Q^{-1}
\ot 1)\phipl \xi \; .
\end{equation}
Combining \eqref{eq.eerstedeel} and \eqref{eq.tweededeel}, we get that the right hand side of \eqref{eq.daargaanwe} equals
$$\sqrt{\frac{\dimq(x+1)}{\dimq(1)\dimq(x)}} \; \; \bigl((Q^{-1} \ot 1)
\phipl \ot \phipr\bigr) (\xi \ot (1 \ot \eta^*)t_x) \; .$$ The formulas for $T^+_r$, $T^-_l$ and $T^-_r$ are proved analogously.
\end{proof}

The proof of Proposition \ref{prop.mainestimate} shall follow immediately from the following two lemmas.

\begin{lemma} \label{lem.firstestimate}
We have the following inequalities for a given $x \geq 1$ and using Notation \ref{not.factintertwin}.
$$(\phipl)^*(Q^{-2}\ot 1)\phipl \; , \quad (\phipr)^*(1 \ot
Q^2)\phipr  \;\; \geq \frac{\dimq(1)\dimq(x)}{\dimq(x+1)} - \frac{\dimq(x-1)}{\dimq(x+1)} \|Q\|^2 \; .$$
\end{lemma}

\begin{proof}
Observe that
$$(\phipl)^*(Q^{-2}\ot 1)\phipl =
\frac{\dimq(1)\dimq(x)}{\dimq(x+1)} \; (t_1^* \ot 1_x)(Q^{-2} \ot p_{x+1})(t_1 \ot 1_x) \; .$$ From a left-handed version of \eqref{eq.wenzl}, we get
$$p_{x+1} = 1 \ot p_x - \frac{\dimq(1)\dimq(x-1)}{\dimq(x)} (1 \ot
p_x)(t_1t_1^* \ot 1_{x-1})(1 \ot p_x) \; .$$ Combining with the previous equality and using the facts
$$t_1^*(Q^{-2} \ot 1)t_1 = 1 \quad\text{and}\quad (t_1^* \ot 1)(Q^{-2}
\ot t_1 t_1^*)(t_1 \ot 1) = \dimq(1)^{-2} Q^{-2} \; ,$$ we obtain the first inequality of the lemma. The second one is proven analogously.
\end{proof}

We prove now a more interesting result, which states that the maps $\sigma$ are \lq far from being intertwiners\rq\ in some sense, at least when
$Q=1$. Observe indeed that the numerical coefficient $(\dimq(x) + 1) / \dimq(x+1)$ in the next statement is always less than $1$.

\begin{lemma} \label{lem.secondestimate}
We have the following inequalities for a given $x \geq 1$ and using Notation \ref{not.factintertwin}.
$$\|(\phipl)^*(Q^{-1-is}\ot 1)\sigma\phipr\| \leq \|Q\| \frac{\dimq(x)
  + 1}{\dimq(x+1)} \; ,$$
for all $s \in \R$.
\end{lemma}

\begin{proof}
First of all, we have
$$(\phipl)^*(Q^{-1-is}\ot 1)\sigma\phipr =
\frac{\dimq(1)\dimq(x)}{\dimq(x+1)} \; (t_1^* \ot 1_x) (Q^{-1-is} \ot p_{x+1}) \si (1_x \ot t_1) \; .$$ From \eqref{eq.recursion}, we get that
$$p_{x+1} = (1 \ot p_x)(p_x \ot 1) -  \frac{(-c)^{x+1}
  \dimq(1)}{\dimq(x)} \; (1 \ot p_x)(t_1 \ot 1_{x-1} \ot t_1^*)(p_x
  \ot 1) \; ,$$
and we easily conclude that
$$(\phipl)^*(Q^{-1-is}\ot 1)\sigma\phipr = \frac{\dimq(x)}{\dimq(x+1)}
\; p_x \bigl( c(1_{x-1} \ot Q^{1+is})\si^* - \frac{(-c)^{x+1}}{\dimq(x+1)} (Q^{-1-is} \ot 1_{x-1})\si \bigr) p_x \; .$$ The lemma follows from this
equality.
\end{proof}

We finally prove Proposition \ref{prop.mainestimate}

\begin{proof}[Proof of Proposition \ref{prop.mainestimate}]
We write $\eta = \sum_x \eta_x$ with $\eta_x \in H_x \ot H_x$ whenever $\eta \in L^2(\cG)$. Take $\eta \in L^2(\cG)$. By Lemmas
\ref{lem.identifications}, \ref{lem.firstestimate} and \ref{lem.secondestimate}, we have, for $x \geq 2$,
$$
\|(T_s^+ \eta)_x \|^2 \geq 2 \bigl( 1 - \|Q\|^2 \bigl( \frac{[x-1]}{[2][x]} + \frac{(1+[x])^2}{[2][x][x+1]} \bigr) \bigr) \|\eta_{x-1}\|^2 = 2 \bigl(
1 - 2 \frac{1+[x]}{[2][x+1]} \|Q\|^2 \bigr) \|\eta_{x-1}\|^2 \; .
$$
We also have
$$\|(T_s^+ \eta)_1 \|^2 = 2(1-|\la(1 \ot Q^{is})t,t\ra|^2)
\|\eta_0\|^2 \; ,$$ and $(T_s^+ \eta)_0 = 0$. By Lemma \ref{lem.identifications}, we have, for all $x \geq 0$,
$$\|(T_s^- \eta)_x \|^2 \leq 4 \|Q\|^2 \frac{[x+1]}{[2][x+2]}
\|\eta_{x+1}\|^2 \; .$$ Suppose now $\frac{\|Q\|}{[2]} \leq \frac{1}{\sqrt{5}}$. We put
$$
C_1 = \sqrt{2} \Bigl(1- 2 \|Q\|^2\frac{1+[2]}{[2][3]} \Bigr)^{1/2} \; , \qquad C_2 = 2 \|Q\| \Bigl(\frac{q}{[2]}\Bigr)^{1/2}
$$
and observe that $C_1$ is well defined and $C_2 < C_1$. So, Proposition \ref{prop.mainestimate} is proved.
\end{proof}

\begin{remark} \label{rem.goedgeteld}
In the proof of Theorem \ref{thm.simplicity}, we need an estimate on the norm of $P$ on $\xi_0^\perp$. With the notations introduced at the end of
the proof of Proposition \ref{prop.mainestimate}, we get
$$\|P\xi\| \leq \Bigl(1-\frac{(C_1 - C_2)^2}{2}\Bigr) \|\xi\|$$
whenever $\xi \in L^2(\cG)$ and $\la \xi_0,\xi \ra = 0$. In order to prove Theorem \ref{thm.simplicity}, we need
\begin{equation*}
\|Q\|^2 \Bigl(1-\frac{(C_1 - C_2)^2}{2}\Bigr) < 1 \; .
\end{equation*}
If $\|Q\|^4 \leq a \Tr(Q)$ with $a \leq 1/\sqrt{5}$, we have
$$\|Q\|^2 \Bigl(1-\frac{(C_1 - C_2)^2}{2}\Bigr) \leq 2a\Bigl(
\frac{1+[2]}{[3]} - q \Bigr) + 2a^{3/4}[2]^{1/4} \sqrt{2q} \sqrt{1-2\frac{1+[2]}{[2][3]}} \; .$$ Taking $a=\frac{3}{8}$ and realizing that $\Tr(Q)
\geq 3$, we conclude that
\begin{equation}\label{eq.daargaathetom}
\|Q\|^2 \|P \xi\| \leq 0,99 \|\xi\| \; ,
\end{equation}
for all $\xi \in \xi_0^\perp$, when $\|Q\|^4 \leq \frac{3}{8} \Tr(Q)$.
\end{remark}

\section{Appendix : approximate commutation of intertwiners} \label{sec.appendix}

In this appendix, we prove several estimates on the representation
theory of $A_o(F)$. Weaker versions of these estimates were proved and
used in \cite{Ver}.

We denote $d_\T(V,W) = \inf\{ \|V - \lambda W \| \mid \lambda \in \T \}$, whenever $V,W$ are in a Banach space.

We fix $F \in \GL(n,\C)$ with $F \overline{F} = c1$ and $c = \pm 1$. We take $0 < q < 1$ such that $\Tr(F^*F) = q + \frac{1}{q}$. We freely use the
explicit choices that can be made for the representation theory of $\cG = A_o(F)$, see Section \ref{sec.factoriality}.

In this appendix, only dealing with the representation theory of $A_o(F)$, all small letters $a,b,c,x,y,z,r,s,\ldots$ denote elements of $\N$, i.e.\
irreducible representations of $A_o(F)$.

\begin{lemma}\label{lem.crucial-improved}
There exists a constant $C > 0$ only depending on $q$ such that for all $a,b,c$ and all $z \in a \ot b$,
\begin{align}
& \| (V(a \ot b,z) \ot 1) p^{z \ot c}_{z+c} - (1 \ot p^{b \ot c}_{b+c})(V(a \ot b,z) \ot 1) \| \leq C q^{(z+b-a)/2} \; , \label{eq.best1} \\
& d_\T\bigl( (V(a \ot b,z) \ot 1)V(z \ot c,z+c) , (1 \ot V(b \ot c,b+c))V(a \ot (b+c),z+c) \bigr) \leq C q^{(z+b-a)/2} \; , \label{eq.best2} \\
& d_\T\bigl( (1 \ot V(b \ot c,b+c)^*)(V(a \ot b,z) \ot 1), V(a \ot (b+c),z+c) V(z \ot c,z+c)^*) \leq C q^{(z+b-a)/2} \; . \label{eq.best3}
\end{align}
\end{lemma}

If we write $z=a+b-2s$, with $0 \leq s \leq \min\{a,b\}$, we have $(z+b-a)/2 = b-s$ and hence, $q^{(z+b-a)/2} \leq q^{-a+b}$.

It is easy to derive \eqref{eq.best2} and \eqref{eq.best3} from \eqref{eq.best1}. Obviously, Lemma \ref{lem.crucial-improved} has a left-handed
analogue.

\begin{lemma} \label{lem.crucial-other-side}
There exists a constant $C > 0$ only depending on $q$ such that for all $a,b,c$ and all $z \in b \ot c$,
\begin{align}
& \| (1 \ot V(b \ot c,z))p^{a \ot z}_{a+z} - (p^{a \ot b}_{a+b} \ot 1)(1 \ot V(b \ot c,z)) \| \leq C q^{(z+b-c)/2} \; , \label{eq.best-a1} \\
& d_\T\bigl( (1 \ot V(b \ot c,z)) V(a \ot z,a+z), (V(a \ot b,a+b) \ot 1) V((a+b) \ot c , a+z) \bigr) \leq C q^{(z+b-c)/2} \; , \label{eq.best-a2} \\
& d_\T\bigl( (V(a \ot b,a+b)^* \ot 1)(1 \ot V(b \ot c,z)) , V((a+b) \ot c , a+z) V(a \ot z,a+z)^* \bigr) \leq C q^{(z+b-c)/2} \; . \label{eq.best-a3}
\end{align}
\end{lemma}

\begin{remark}
It is possible to prove Lemma \ref{lem.crucial-improved} from explicit formulae for the quantum $6j$-symbols of $\SU_{q}(2)$. We give a more direct
approach, for which we only need to know the quantum $3j$-symbols (i.e.\ the coefficient appearing in \eqref{eq.innerproduct}).
\end{remark}

Before giving the proof of Lemma \ref{lem.crucial-improved}, we introduce several notations and lemmas. It follows from \cite{BDRV} (using
\cite{banica1,banica2}) that $A_o(F)$ is monoidally equivalent with $\SU_{-cq}(2)$, where $q$ is as before and $F \overline{F} = c 1$, $c = \pm 1$.
So, we can perform all computations on the intertwiners as if we are dealing with the representation theory of $\SU_{-cq}(2)$.

Recall the explicit choices that can be made for the representation theory of $A_o(F)$ in Section \ref{sec.factoriality}.

\begin{lemma} \label{lem.approxcommute}
There exists a constant $C$, only depending on $q$, such that
$$\| (p_{a+b} \ot 1_c)(1_a \ot p_{b+c}) - p_{a+b+c} \| \leq C
q^{b} \; .$$
\end{lemma}
\begin{proof}
From \eqref{eq.recursion}, it follows that
\begin{equation}\label{eq.one}
p_{b+c+1} = (p_b \ot p_{c+1})p_{b+c+1} = \Bigl(1_b \ot p_{c+1}
- \frac{[2][b]}{[b+c+1]} (p_b \ot p_{c+1})(1_{b-1} \ot t \ot
1_c \ot t^*) \Bigr) (p_{b+c} \ot 1) \; .
\end{equation}
In the same way, it follows that
\begin{align*}
p_{a+b+c+1} &= (p_a \ot p_b \ot p_{c+1}) p_{a+b+c+1}  \\ & =\Bigl(1_{a+b}
  \ot p_{c+1} - \frac{[2][a]}{[a+b+c+1]} (p_a \ot p_b  \ot
  p_{c+1})(1_{a-1} \ot t \ot 1_{b+c} \ot t^*)
\\ &\hspace{2.5cm} - \frac{[2][a+b]}{[a+b+c+1]} (p_a \ot p_b  \ot
  p_{c+1})(1_{a+b-1} \ot t \ot 1_c \ot
  t^*)\Bigr)(p_{a+b+c} \ot 1) \; .
\end{align*}
Since both $\frac{[2][a]}{[a+b+c+1]}$ and the difference $\bigl| \frac{[2][a+b]}{[a+b+c+1]} - \frac{[2][b]}{[b+c+1]}\bigr| =
\frac{[2][a][c+1]}{[b+c+1][a+b+c+1]}$ can be estimated by $C q^{b+c}$ for a constant $C$ only depending on $q$, we find that
\begin{equation}\label{eq.two}
p_{a+b+c+1} \approx \Bigl( 1_{a+b} \ot p_{c+1} -
\frac{[2][b]}{[b+c+1]} (1_a \ot p_b \ot p_{c+1})(1_{a+b-1} \ot t \ot
1_c \ot t^*)\Bigr)(p_{a+b+c} \ot 1)
\end{equation}
with error $\leq C q^{b+c}$.

Put now $\eps(a,b,c) = \|(1_a \ot p_{b+c})(p_{a+b} \ot 1) - p_{a+b+c}
\|$. Using \eqref{eq.one}, we find that
\begin{align*}
(1_a & \ot p_{b+c+1})(p_{a+b} \ot 1_{c+1}) \\ & = \Bigl( 1_{a+b} \ot p_{c+1} -
\frac{[2][b]}{[b+c+1]} (1_a \ot p_b \ot p_{c+1})(1_{a+b-1} \ot t \ot
1_c \ot t^*)\Bigr) (1_a \ot p_{b+c} \ot 1)(p_{a+b} \ot 1_{c+1}) \\ &
\approx \Bigl( 1_{a+b} \ot p_{c+1} -
\frac{[2][b]}{[b+c+1]} (1_a \ot p_b \ot p_{c+1})(1_{a+b-1} \ot t \ot
1_c \ot t^*)\Bigr) (p_{a+b+c} \ot 1)
\end{align*}
with error $\leq \eps(a,b,c)(1 + D q^{c})$, because we can find $D$ such that $\frac{[2][b]}{[b+c+1]} \leq D q^{c}$. Combining with \eqref{eq.two},
we find that
$$\eps(a,b,c+1) \leq \eps(a,b,c)(1 + Dq^{c}) + Cq^{b+c} \; .$$
One concludes by induction that
$$\eps(a,b,c) \leq \sum_{k=0}^{c-1} \Bigl(C q^{b+k} \prod_{j=k+1}^{c-1} (1 + Dq^j) \Bigr) \; .$$
It follows that
$$\eps(a,b,c) \leq q^{b} \Bigl(\prod_{k=0}^\infty
(1+Dq^{k})\Bigr) \Bigl(\sum_{k=0}^\infty C q^{k} \Bigr) \; .$$ This concludes the proof of the lemma.
\end{proof}

\begin{lemma}
There exists a constant $C$ only depending on $q$ such that
$$1 \leq \frac{\qbinom{a+r}{r}
  \qbinom{r+b}{r}}{\qbinom{a+b+r}{r}} \leq C$$
for all $a,b,r$.
\end{lemma}
\begin{proof}
It is easy to find a constant $C$ such that for all $a,b$ and all $k
\geq 1$, we have
\begin{equation}\label{eq.interessant}
1 \leq \frac{[a+k] [b+k]}{[a+b+k] [k]} \leq 1 + C q^{2k} \; .
\end{equation}
Taking the product for $k$ running from $1$ to $r$, we obtain the result because $\prod_{k=0}^\infty (1 + C q^{k}) < +\infty$.
\end{proof}

\begin{proof}[Proof of Lemma \ref{lem.crucial-improved}]
It suffices to prove \eqref{eq.best1}. We introduce the notation
\begin{equation}\label{eq.notation}
C(a,b,r) = \frac{[r+1] \qbinom{a+r}{r}
  \qbinom{r+b}{r}}{\qbinom{a+b+r+1}{r}} \; .
\end{equation}
We identify
\begin{align}
& (V((a+s) \ot (s + b),a+b)^* \ot 1) (1_{a+s}  \ot p^{(s+b) \ot c}_{s+b+c}) (V((a+s) \ot (s + b),a+b) \ot 1) \label{eq.tussenres} \\ & = C(a,b,s)
(p_{a+b} \ot 1_c)(1_a \ot p_{b+c}) \;  (1_a \ot t^*_s \ot 1_{b+c})(p_{a+s} \ot p_{s+b+c}) (1_a \ot t_s \ot 1_{b+c}) \; (1_a \ot p_{b+c}) (p_{a+b} \ot
1_c) \; . \notag
\end{align}
But,
$$
C(a,b,s) (1_a \ot t^*_s \ot 1_{b+c})(p_{a+s} \ot p_{s+b+c}) (1_a \ot t_s \ot 1_{b+c}) = \sum_{z \in a \ot (b+c)} \lambda_z p^{a \ot (b+c)}_z \; .
$$
From \eqref{eq.innerproduct}, we get that $\lambda_{a+b+c} = C(a,b,s) / C(a,b+c,s)$. Since $C(a,b,s)$ is uniformly bounded from above, we get a
constant $D$ only depending on $q$ such that $\lambda_z \leq D$ for all $z$. Using Lemma \ref{lem.approxcommute}, we find a constant $E$ such that,
for all $z < a+b+c$,
$$\| p^{a \ot {b+c}}_z \; (1_a \ot p_{b+c}) (p_{a+b} \ot 1_c) \| \leq E q^{b} \; .$$
As in \eqref{eq.interessant}, we find that
$$1-q^{2(b+k)} \leq \frac{[b+k] [a+b+c+1+k]}{[a+b+1+k][b+c+k]} \leq
1 \; .$$ Taking the product for $k$ running from $1$ to $s$, we find a constant $G$ such that
$$1-Gq^{2b} \leq \frac{C(a,b,s)}{C(a,b+c,s)} \leq 1 \; .$$
Combining all these estimates with \eqref{eq.tussenres}, we have shown the existence of a constant $C$, only depending on $q$, such that
$$\bigl\|
(V((a+s) \ot (s + b),a+b)^* \ot 1) (1_{a+s}  \ot p^{(s+b) \ot c}_{s+b+c}) (V((a+s) \ot (s + b),a+b) \ot 1) - p_{a+b+c} \bigr\| \leq (C q^{b})^2 \;
.$$ It follows that
$$\bigl\| (1_{a+s}  \ot p^{(s+b) \ot c}_{s+b+c}) (V((a+s) \ot (s + b),a+b) \ot 1) - (V((a+s) \ot (s + b),a+b) \ot 1) p_{a+b+c}\bigr\| \leq C\sqrt{2}
q^{b} \; .$$ This is the formula we had to prove.
\end{proof}

\begin{lemma} \label{lem.encore-une}
There exists a constant $D>0$ only depending on $q$ such that
$$D \|\xi\| \leq \|(p_{a+s} \ot p_{s+b})(1_a \ot t_s \ot 1_b) \xi \| \leq \|\xi\|
$$
for all $a,b,s$ and $\xi \in H_a \ot H_b$.
\end{lemma}
\begin{proof}
Use the notation \eqref{eq.notation}. Since
$$C(a,b,r) = \frac{\qbinom{a+b+r}{r} \qbinom{1+r}{r}}{\qbinom{a+b+r+1}{r}} \frac{\qbinom{a+r}{r}
  \qbinom{r+b}{r}}{\qbinom{a+b+r}{r}} \; ,$$
it follows from the previous lemma that there exist constants $C_1$ and $C_2$ such that $C_1 \leq C(a,b,r) \leq C_2$. Observe
\begin{align*}
(p_{a+x+s} \ot p_{s+x+b}) & (1_{a+x} \ot t_s \ot 1_{x+b}) V((a+x) \ot (x+b),a+b) \\ &= C(a,b,x)^{1/2}
(p_{a+x+s} \ot p_{s+x+b})(1_{a+x} \ot t_s \ot 1_{x+b})(1_a \ot t_x \ot 1_b)p_{a+b} \\
&= D(x,s)^{-1/2} C(a,b,x)^{1/2} C(a,b,x+s)^{-1/2} V((a+x+s) \ot (s+x+b),a+b) \;,
\end{align*}
where $D(x,s) = [x+1] [s+1] [x+s+1]^{-1}$. Since also $D(x,s)$ lies between two constants, the lemma is proved.
\end{proof}

\end{document}